# Campbell equilibrium equation and pseudo-likelihood estimation for non-hereditary Gibbs point processes

DAVID DEREUDRE[1] and FRÉDÉRIC LAVANCIER[2]

[1]*LAMAV, Université de Valenciennes et du Hainaut-Cambrésis, Le Mont Houy, 59313 Valenciennes Cedex 09, France. E-mail: david.dereudre@univ-valenciennes.fr*
[2]*Université de Nantes, Laboratoire de Mathématiques Jean Leray, Unité Mixte de Recherche CNRS 6629, UFR Sciences et Techniques, 2 rue de la Houssinière – BP 92208 – F-44322 Nantes Cedex, France. E-mail: frederic.lavancier@univ-nantes.fr*

In this paper, we study Gibbs point processes involving a hardcore interaction which is not necessarily hereditary. We first extend the famous Campbell equilibrium equation, initially proposed by Nguyen and Zessin [*Math. Nachr.* **88** (1979) 105–115], to the non-hereditary setting and consequently introduce the new concept of removable points. A modified version of the pseudo-likelihood estimator is then proposed, which involves these removable points. We consider the following two-step estimation procedure: first estimate the hardcore parameter, then estimate the smooth interaction parameter by pseudo-likelihood, where the hardcore parameter estimator is plugged in. We prove the consistency of this procedure in both the hereditary and non-hereditary settings.

## 1. Introduction

Gibbs point processes first appeared in statistical physics for the description of large interacting particle systems. They are now used in many other fields such as biology, medicine, agronomy and astronomy. Gibbs models rely on interaction potential functions that must be properly estimated in applications. We are particularly interested in the case where the potential function involves a hardcore interaction. This means that some point configurations are forbidden with respect to this potential function or, equivalently, that some point configurations have infinite energy. This constraint appears in many classical models of Gibbs measures (for example, the hard ball model).

To the best of our knowledge, in all previous studies involving a hardcore interaction, the interaction is assumed to be hereditary (see [15]). As recalled in Definition 1, an interaction is hereditary if, for every forbidden point configuration $\gamma$ in $\mathbb{R}^d$ and every $x$ in $\mathbb{R}^d$, the configuration $\gamma + \delta_x$ (that is, the configuration composed of the union of the point $x$ with the points of $\gamma$) remains forbidden. Note that any interaction involving no hardcore part is necessarily hereditary.

In this paper, we aim to study Gibbs measures in the presence of a hardcore interaction that is not necessarily hereditary. Indeed, it seems natural to encounter some non-hereditary interactions, as Examples 5.2 and 5.3 in Section 5 show. The first example is concerned with random fields of geometric objects subjected to a hardcore interaction. Such a model is developed in [5] with the study of random tessellations with geometric hardcore interaction. In this model, each cell of the tessellation is forced to have a radius less than a fixed constant $\alpha$. It is clear that this kind of interaction is not hereditary. In fact, when one adds a new point inside a forbidden large cell, the new tessellation may become authorized. The second example is a forced clustering process. In this model, each point is forced to have at least $k$ neighbors at a distance less than $\alpha$. Here, again, this interaction is clearly not hereditary.

In the first part of the paper, we introduce some notation, recall the definition of Gibbs measures on $\mathbb{R}^d$ and prove some preliminary results about non-hereditary Gibbs processes. The main problem, in the non-hereditary case, is that the energy of a point $x$ in a configuration $\gamma$ is not always defined. Indeed, the energy of $\gamma - \delta_x$ may be locally infinite, even if the energy of $\gamma$ is locally finite. In this case, the energy of $x$ in $\gamma$ would be negative infinity, which makes no sense. We must then introduce the set $\mathcal{R}(\gamma)$ of removable points in $\gamma$ (see Definition 2): $x$ in $\gamma$ is said to be *removable* if the energy of $\gamma - \delta_x$ is locally finite.

The second part of the paper is dedicated to the elaboration of a Campbell equilibrium equation for non-hereditary Gibbs measures. The classical Campbell equilibrium equation, proposed by Nguyen and Zessin [17], is only valid in the hereditary case. This equation is fundamental for the understanding of Gibbs measures. It also provides major statistical applications: the Takacs–Fiksel and pseudo-likelihood estimation procedures rely on the Nguyen–Zessin equation (see, for instance, [15] for an exposition). Moreover, the recent definition of residuals for spatial point processes arises from this equation (see [3]). The Introduction of the concept of removable points allows us to extend the Campbell equilibrium equation when the heredity of the interaction is not assumed. Let $\mu$ be a Gibbs measure and $\mathcal{C}_\mu^!$ its reduced Campbell measure. We prove in Proposition 2 that

$$\mathbb{1}_{x \in \mathcal{R}(\gamma + \delta_x)} \mathcal{C}_\mu^!(\mathrm{d}x, \mathrm{d}\gamma) = \mathrm{e}^{-h(x,\gamma)} \lambda \otimes \mu(\mathrm{d}x, \mathrm{d}\gamma), \qquad (1)$$

where $h$ denotes the local energy of $x$ in $\gamma$ and $\lambda$ is the Lebesgue measure on $\mathbb{R}^d$.

Next, we consider the estimation problem. Several statistical procedures have been proposed to estimate the potential function of a hereditary Gibbs measure. Among them, the common parametric methods are the maximum likelihood estimation and the pseudo-likelihood procedure. In addition, the Takacs–Fiksel method, some Bayesian procedures, and semi- and non-parametric methods have been developed. We refer to [15] for a recent review of these estimation procedures.

In Section 4, assuming the interaction is parametric, we focus on the pseudo-likelihood method. The direct parametric method consists of the maximum likelihood procedure, but it requires the estimation, by simulations, of an unknown normalizing constant. Recent advances, including Monte Carlo and perfect simulation, deal with these simulations efficiently (see [4, 7, 13]). Yet, from a practical point of view, maximum likelihood remains computationally much more intensive than pseudo-likelihood. Therefore, pseudo-likelihood estimation is often used as a first step before maximum likelihood estimation. For the same reason, it may be preferred over maximum likelihood estimation when fitting several alternative models to a data set. Moreover, in the stationary case, some asymptotic results exist for the maximum likelihood procedure (see the review in [15]), but they rely on rather restrictive assumptions that seem difficult to avoid. The pseudo-likelihood method has the advantage of satisfying the usual asymptotic theory for a large class of interactions; see [2, 8–10, 14]. In these articles, under very general assumptions, the consistency of the pseudo-likelihood method is proved, as well as the asymptotic normality of the induced estimators. In the same spirit, in Section 4, we prove the consistency of the pseudo-likelihood procedure when the interaction is not necessarily hereditary.

More precisely, we suppose that the interaction depends on two parameters, $\alpha$ and $\theta$. The parameter $\alpha$ is devoted to the hardcore interaction, whereas $\theta$ parameterizes the classical smooth interaction. In our non-hereditary setting, the pseudo-likelihood contrast function is defined as

$$PLL_{\Lambda_n}(\gamma, \alpha, \theta) = \frac{1}{|\Lambda_n|} \left[ \int_{\Lambda_n} \exp(-h^{\alpha,\theta}(x,\gamma)) \, \mathrm{d}x + \sum_{x \in \mathcal{R}^{\alpha}(\gamma) \cap \Lambda_n} h^{\alpha,\theta}(x, \gamma - \delta_x) \right], \qquad (2)$$

where $\Lambda_n$ denotes the observation domain of the sample. This pseudo-likelihood contrast function is the classical one found in previous studies, but in order to remain well defined in the non-hereditary setting, the sum is restricted to the removable points $\mathcal{R}^{\alpha}(\gamma)$ (this set depends only on the hardcore parameter $\alpha$, as implied by assumption S1 in Section 4). We consider a classical two-step estimation procedure. We first estimate $\alpha$ in a consistent way (see Proposition 3). We then estimate $\theta$ by maximizing $PLL_{\Lambda_n}$, where $\alpha$ is replaced by its estimator. Theorem 2 establishes the consistency of both estimators of $\alpha$ and $\theta$ resulting from this procedure. Let us note that even if the interaction is hereditary, to the best of our knowledge, no consistency proof has previously been supplied for this classical two-step procedure, although it has been widely used in practice. This is, for instance, the case for the classical hardball model that we present in Section 5.

Section 4 is devoted to the consistency of the estimators, but their asymptotic normality is not addressed. However, there are no major difficulties in obtaining the asymptotic normality of the estimator of $\theta$ when $\alpha$ is known, provided the interactions are finite range. The tools for proving this are essentially the same as in [10] and, more recently, in [2]. In the situation where the hardcore parameter $\alpha$ is estimated, the asymptotic normality of the estimator of $\theta$ seems more difficult to obtain.

Section 5 contains some examples of parametric interactions involving a hardcore part. The first is the hardball model with pairwise step interaction. It is a hereditary model. The second is the hardcore Delaunay tessellations model, which is not hereditary. The third example is a forced clustering process with interaction, which is also not hereditary.

For these examples, we prove that the assumptions of Theorem 2 are satisfied, leading to the consistency of the joint estimation of their parameters.

## 2. Definitions and notation

### 2.1. State spaces and reference measures

Let $d$ be a fixed integer greater than 1. $\mathbb{R}^d$ denotes the $d$-dimensional Euclidean real space endowed with the Borel $\sigma$-algebra $\sigma(\mathbb{R}^d)$. $\mathcal{B}(\mathbb{R}^d)$ is the set of bounded Borel sets on $\mathbb{R}^d$. $\mathcal{M}(\mathbb{R}^d)$ is the set of simple integer-valued measures $\gamma$ on $\mathbb{R}^d$, that is, for every $\Lambda \in \mathcal{B}(\mathbb{R}^d)$, $\gamma(\Lambda) \in \mathbb{N}$ and for every $x \in \mathbb{R}^d$, $\gamma(\{x\}) \leq 1$. $\mathcal{M}(\mathbb{R}^d)$ is endowed with the $\sigma$-algebra $\sigma(\mathcal{M}(\mathbb{R}^d))$ generated by the sets $\{\gamma \in \mathcal{M}(\mathbb{R}^d), \gamma(\Lambda) = n\}$, $n \in \mathbb{N}$, $\Lambda \in \mathcal{B}(\mathbb{R}^d)$. So, $\mathcal{P}(\mathcal{M}(\mathbb{R}^d))$ denotes the set of probability measures on $\mathcal{M}(\mathbb{R}^d)$ for the $\sigma$-algebra $\sigma(\mathcal{M}(\mathbb{R}^d))$. Any measure $\gamma \in \mathcal{M}(\mathbb{R}^d)$ has the following representation:

$$\gamma = \sum_{i \in \mathcal{I}} \delta_{x_i},$$

where $\mathcal{I}$ is a subset of $\mathbb{N}$, $(x_i)_{i \in \mathcal{I}}$ are elements of $\mathbb{R}^d$ and $\delta_x$ is the Dirac measure at $x$. We write $x \in \gamma$ if $\gamma(\{x\}) > 0$.

Letting $\gamma$ be in $\mathcal{M}(\mathbb{R}^d)$ and $\Lambda$ be a Borel set in $\mathbb{R}^d$, we denote by $\gamma_\Lambda$ the projection of $\gamma$ on $\Lambda$, which is just the measure $\sum_{x \in \gamma \cap \Lambda} \delta_x$.

Denoting by $\lambda$ the Lebesgue measure on $\mathbb{R}^d$, $\pi$ stands for the Poisson process on $\mathbb{R}^d$ with intensity $\lambda$. It is a probability measure on $\mathcal{M}(\mathbb{R}^d)$. For every $\Lambda \in \mathcal{B}(\mathbb{R}^d)$, $\pi_\Lambda$ denotes the Poisson process on $\Lambda$ with intensity $\mathbb{1}_\Lambda \lambda$. It is a probability measure on $\mathcal{M}(\Lambda)$.

### 2.2. Interaction

We define the interaction energy in a general setting, as in [12]. We consider the most general specifications of Gibbs kernels with local densities having an exponential form. We do not assume that the local densities come from a multibody interaction potential. This general point of view allows us to deal with the non-hereditary case (see Definition 1 below).

A set of measurable functions $(H_\Lambda)_{\Lambda \in \mathcal{B}(\mathbb{R}^d)}$ from $\mathcal{M}(\mathbb{R}^d)$ to $\mathbb{R} \cup \{+\infty\}$ defines a family of energies if, for every $\Lambda \subset \Lambda'$ in $\mathcal{B}(\mathbb{R}^d)$, there exists a measurable function $\varphi_{\Lambda,\Lambda'}$ from $\mathcal{M}(\mathbb{R}^d)$ to $\mathbb{R} \cup \{+\infty\}$ such that

$$\forall \gamma \in \mathcal{M}(\mathbb{R}^d) \qquad H_{\Lambda'}(\gamma) = H_\Lambda(\gamma) + \varphi_{\Lambda,\Lambda'}(\gamma_{\Lambda^c}). \tag{3}$$

(3) is equivalent to (6.11) and (6.12) from page 92 of [12]. In particular, we have the following property:

$$\forall \Lambda \subset \Lambda' \text{ in } \mathcal{B}(\mathbb{R}^d) \text{ and } \forall \gamma \in \mathcal{M}(\mathbb{R}^d) \qquad H_\Lambda(\gamma) = +\infty \quad \Rightarrow \quad H_{\Lambda'}(\gamma) = +\infty. \tag{4}$$

Physically, $H_\Lambda(\gamma)$ is the energy of $\gamma_\Lambda$ inside $\Lambda$ given the configuration $\gamma_{\Lambda^c}$ outside $\Lambda$. Let us now discuss the problem concerning heredity.

**Definition 1.** *A family of energies is* hereditary *if*

$$\forall \Lambda \in \mathcal{B}(\mathbb{R}^d), \forall \gamma \in \mathcal{M}(\mathbb{R}^d) \text{ and } \forall x \in \Lambda \qquad H_\Lambda(\gamma) = +\infty \quad \Rightarrow \quad H_\Lambda(\gamma + \delta_x) = +\infty. \qquad (5)$$

This assumption (5) is often needed in many papers (see, for example, [16, 17]). Let us point out that the non-heredity of a family of energies necessarily comes from a hardcore part (that is, configurations with infinite energy). From now on, we will never assume that the energy is hereditary. When we invoke the "non-hereditary case", it means the general setting of hereditary or non-hereditary energies.

Let us introduce the new concept of removable points, which is fundamental in the non-hereditary case.

**Definition 2.** *Let $\gamma$ be in $\mathcal{M}(\mathbb{R}^d)$ and $x$ be a point of $\gamma$. $x$ is* removable *from $\gamma$ if*

$$\exists \Lambda \in \mathcal{B}(\mathbb{R}^d) \qquad \text{such that } x \in \Lambda \quad \text{and} \quad H_\Lambda(\gamma - \delta_x) < +\infty. \qquad (6)$$

*We denote by $\mathbb{R}(\gamma)$ the set of removable points in $\gamma$.*

When the configuration $\gamma$ has a locally finite energy, this definition can be simplified, as proved in Proposition 1. A configuration $\gamma$ in $\mathcal{M}(\mathbb{R}^d)$ has a locally finite energy if for every $\Lambda \in \mathcal{B}(\mathbb{R}^d)$, the energy $H_\Lambda(\gamma)$ is finite. We denote by $\mathcal{M}_\infty(\mathbb{R}^d)$ the space of configurations which have a locally finite energy.

**Proposition 1.** *If $\gamma$ is in $\mathcal{M}_\infty(\mathbb{R}^d)$ and $x$ is a point of $\gamma$, then $x$ is removable from $\gamma$ if and only if $\gamma - \delta_x$ is in $\mathcal{M}_\infty(\mathbb{R}^d)$.*

**Proof.** Let $\gamma$ be in $\mathcal{M}_\infty(\mathbb{R}^d)$ and $x$ be a removable point of $\gamma$. There exists $\Lambda \in \mathcal{B}(\mathbb{R}^d)$ such that $x$ is in $\Lambda$ and $H_\Lambda(\gamma - \delta_x)$ is finite. Let us show that $H_{\Lambda'}(\gamma - \delta_x) < +\infty$ for every bounded set $\Lambda'$ in $\mathbb{R}^d$. Thanks to (3), we have

$$H_{\Lambda' \cup \Lambda}(\gamma - \delta_x) = H_\Lambda(\gamma - \delta_x) + \varphi_{\Lambda, \Lambda \cup \Lambda'}((\gamma - \delta_x)_{\Lambda^c})$$
$$= H_\Lambda(\gamma - \delta_x) + \varphi_{\Lambda, \Lambda \cup \Lambda'}(\gamma_{\Lambda^c}).$$

Since $\gamma$ is in $\mathcal{M}_\infty(\mathbb{R}^d)$, $\varphi_{\Lambda, \Lambda \cup \Lambda'}(\gamma_{\Lambda^c})$ is finite. Thus, $H_{\Lambda' \cup \Lambda}(\gamma - \delta_x)$ is finite and, thanks to (4), $H_{\Lambda'}(\gamma - \delta_x)$ is also finite. Therefore, $\gamma - \delta_x$ is in $\mathcal{M}_\infty(\mathbb{R}^d)$.

The converse is obvious. $\square$

We can now define the local energy of a removable point $x$ in a configuration $\gamma$.

**Definition 3.** *Let $x$ be a removable point in a configuration $\gamma$ in $\mathcal{M}(\mathbb{R}^d)$. The local energy of $x$ in $\gamma - \delta_x$ is defined as*

$$h(x, \gamma - \delta_x) = H_\Lambda(\gamma) - H_\Lambda(\gamma - \delta_x), \qquad (7)$$

*where $\Lambda$ is a bounded set containing $x$ such that $H_\Lambda(\gamma - \delta_x)$ is finite.*

Notice that this definition is valid and does not depend on the choice of $\Lambda$. Indeed, according to Definition 2, at least one such $\Lambda$ exists. Besides, let us suppose that there is another $\Lambda'$ containing $x$ and such that $H_{\Lambda'}(\gamma - \delta_x)$ is finite. Writing $\Lambda'' = \Lambda \cap \Lambda'$, we have, from (3),

$$\begin{aligned} H_\Lambda(\gamma) - H_\Lambda(\gamma - \delta_x) &= H_{\Lambda''}(\gamma) + \varphi_{\Lambda'',\Lambda}(\gamma_{\Lambda''^c}) - H_{\Lambda''}(\gamma - \delta_x) - \varphi_{\Lambda'',\Lambda}((\gamma - \delta_x)_{\Lambda''^c}) \\ &= H_{\Lambda''}(\gamma) + \varphi_{\Lambda'',\Lambda}(\gamma_{\Lambda''^c}) - H_{\Lambda''}(\gamma - \delta_x) - \varphi_{\Lambda'',\Lambda}(\gamma_{\Lambda''^c}) \\ &= H_{\Lambda''}(\gamma) - H_{\Lambda''}(\gamma - \delta_x), \end{aligned}$$

which is equal, thanks to the same calculations, to $H_{\Lambda'}(\gamma) - H_{\Lambda'}(\gamma - \delta_x)$. Therefore, $h(x, \gamma - \delta_x)$ is well defined and belongs to $\mathbb{R} \cup \{\infty\}$.

Finally, for every $\gamma$ in $\mathcal{M}_\infty(\mathbb{R}^d)$ and every $x$ in $\mathbb{R}^d$, the local energy of $x$ in $\gamma$ is always defined since $h(x, \gamma) = h(x, \gamma + \delta_x - \delta_x)$ and $H_\Lambda(\gamma)$ is finite for all bounded Borel sets $\Lambda$ in $\mathbb{R}^d$.

## 2.3. Gibbs states

In this subsection, we are in a position to define the Gibbs states. We need to introduce the notion of specifications, as in [12]. Let us make an integrability assumption about the family of energies, which is equivalent to (6.8) in [12]. We say that the family of energies $(H_\Lambda)$ is *integrable* if, for every $\Lambda$ in $\mathcal{B}(\mathbb{R}^d)$ and every $\gamma$ in $\mathcal{M}_\infty(\mathbb{R}^d)$, we have

$$0 < \int_{\mathcal{M}(\Lambda)} e^{-H_\Lambda(\gamma'_\Lambda + \gamma_{\Lambda^c})} \pi_\Lambda(d\gamma'_\Lambda) < +\infty. \tag{8}$$

The second inequality in (8) is, in general, guaranteed by the stability of the potential. We will assume this stability in H3 below. The first inequality is obvious in the classical hereditary setting. In the non-hereditary setting, it remains true under reasonable assumptions (see, for instance, [5]).

Under this integrability assumption, we are able to define the kernels for the Gibbs structure (see (6.6) in [12]). For every $\Lambda$ in $\mathcal{B}(\mathbb{R}^d)$, the kernel $\Xi_\Lambda$ on $\mathcal{P}(\mathcal{M}_\infty(\mathbb{R}^d)) \times \mathcal{M}_\infty(\mathbb{R}^d)$ is defined by

$$\begin{aligned} \Xi_\Lambda(f, \gamma) &= \int_{\mathcal{M}(\mathbb{R}^d)} f(\gamma') \Xi_\Lambda(d\gamma', \gamma) \\ &:= \int_{\mathcal{M}(\Lambda)} f(\gamma'_\Lambda + \gamma_{\Lambda^c}) \frac{1}{Z_\Lambda(\gamma_{\Lambda^c})} e^{-H_\Lambda(\gamma'_\Lambda + \gamma_{\Lambda^c})} \pi_\Lambda(d\gamma'_\Lambda), \end{aligned} \tag{9}$$

where $f$ is just a bounded measurable test function and $Z_\Lambda(\gamma_{\Lambda^c})$ is the normalization constant defined by

$$Z_\Lambda(\gamma_{\Lambda^c}) = \int_{\mathcal{M}(\Lambda)} e^{-H_\Lambda(\gamma'_\Lambda + \gamma_{\Lambda^c})} \pi_\Lambda(d\gamma'_\Lambda).$$

Note that from (8), $0 < Z_\Lambda(\gamma_{\Lambda^c}) < +\infty$ and therefore the kernels are well defined. Moreover, due to (3), they are compatible in the sense of (2.14) in [12], which means that for every $\Lambda \subset \Lambda'$, every $\gamma$ in $\mathcal{M}_\infty(\mathbb{R}^d)$ and every bounded measurable function $f$,

$$\int_{\mathcal{M}(\mathbb{R}^d)} f(\gamma') \Xi_{\Lambda'}(\mathrm{d}\gamma', \gamma) = \int_{\mathcal{M}(\mathbb{R}^d)} \int_{\mathcal{M}(\mathbb{R}^d)} f(\gamma'') \Xi_\Lambda(\mathrm{d}\gamma'', \gamma') \Xi_{\Lambda'}(\mathrm{d}\gamma', \gamma). \qquad (10)$$

Now let us give the definition of Gibbs measures (see (2.15) in [12]).

**Definition 4.** *A probability measure $\mu$ on $\mathcal{M}(\mathbb{R}^d)$ is a Gibbs measure for the family of integrable energies $(H_\Lambda)$ if, for every $\Lambda$ in $\mathcal{B}(\mathbb{R}^d)$ and every bounded measurable function $f$ from $\mathcal{M}(\mathcal{R}^d)$ to $\mathbb{R}$, we have*

$$\int_{\mathcal{M}(\mathbb{R}^d)} f(\gamma) \mu(\mathrm{d}\gamma) = \int_{\mathcal{M}(\mathbb{R}^d)} \int_{\mathcal{M}(\mathcal{R}^d)} f(\gamma') \Xi_\Lambda(\mathrm{d}\gamma', \gamma) \mu(\mathrm{d}\gamma). \qquad (11)$$

We denote by $\mathcal{G}$ the set of stationary Gibbs measures. Equations of the form (11) are called DLR equations, where DLR stands for Dobrushin, Landford and Ruelle. They may be rewritten in the following way: for $\mu$-a.e. $\gamma$ and for every bounded set $\Lambda$ in $\mathcal{B}(\mathbb{R}^d)$,

$$\mu(\cdot | \gamma_{\Lambda^c}) = \Xi_\Lambda(\cdot, \gamma).$$

**Remark 1.** From (11), we deduce that the support of $\mu$ is included in $\mathcal{M}_\infty(\mathbb{R}^d)$.

## 3. Campbell equilibrium equation for non-hereditary Gibbs point processes

In this section, we develop a Campbell equilibrium equation for non-hereditary Gibbs point processes. In [17], the authors give an equation using the reduced Campbell measure to characterize hereditary Gibbs point processes. This Nguyen–Zessin equation (12) is very well known and is used in many works concerning Gibbs processes:

$$\mathcal{C}_\mu^! = \mathrm{e}^{-h} \lambda \otimes \mu. \qquad (12)$$

Unfortunately, this formula is not valid in the non-hereditary case, as explained in Remark 2 below. Therefore, we propose to generalize it, in Proposition 2, to the non-hereditary case. In the hereditary setting, it simply becomes the classical Nguyen–Zessin equation (12). The concept of removable points introduced earlier is the key for this generalization.

First, let us introduce the definition of the reduced Campbell measure (see, for example, [11], page 225).

**Definition 5.** *Let $\mu$ be a probability measure on $\mathcal{M}(\mathbb{R}^d)$. The reduced Campbell measure $\mathcal{C}^!_\mu$ on $\mathbb{R}^d \times \mathcal{M}(\mathbb{R}^d)$ is defined by*

$$\mathcal{C}^!_\mu(f) = \int_{\mathcal{M}(\mathbb{R}^d)} \int_{\mathbb{R}^d} f(x, \gamma - \delta_x) \gamma(\mathrm{d}x) \mu(\mathrm{d}\gamma),$$

*where $f$ is a bounded non-negative measurable function from $\mathbb{R}^d \times \mathcal{M}(\mathbb{R}^d)$ to $\mathbb{R}$.*

**Remark 2.** We shall explain why the Nguyen–Zessin equation is not valid in general for a non-hereditary Gibbs point process. Let us show that $\mathcal{C}^!_\mu$ is not absolutely continuous with respect to $\lambda \otimes \mu$, which implies that (12) is not satisfied. For this, since the support of $\mu$ is $\mathcal{M}_\infty(\mathbb{R}^d)$, it is sufficient to prove that the support of $\mathcal{C}^!_\mu$ is not included in $\mathbb{R}^d \times \mathcal{M}_\infty(\mathbb{R}^d)$. Indeed, from Definition 5,

$$\mathcal{C}^!_\mu([\mathbb{R}^d \times \mathcal{M}_\infty(\mathbb{R}^d)]^c) = \int_{\mathcal{M}(\mathbb{R}^d)} \int_{\mathbb{R}^d} \mathbb{1}_{\gamma - \delta_x \notin \mathcal{M}_\infty(\mathbb{R}^d)} \gamma(\mathrm{d}x) \mu(\mathrm{d}\gamma),$$

where $[A]^c$ denotes the complementary set of $A$ in $\mathbb{R}^d \times \mathcal{M}(\mathbb{R}^d)$. But, for a non-removable point $x$ in $\gamma$, $\gamma - \delta_x$ is not in $\mathcal{M}_\infty(\mathbb{R}^d)$ (see Proposition 1). Since non-hereditary Gibbs point processes contain some non-removable points, $\mathcal{C}^!_\mu([\mathbb{R}^d \times \mathcal{M}_\infty(\mathbb{R}^d)]^c)$ does not vanish in general and the support of $\mathcal{C}^!_\mu$ is not included in $\mathbb{R}^d \times \mathcal{M}_\infty(\mathbb{R}^d)$. This is exactly the situation in Examples 2 and 3 in the last section.

Let us now present our generalization for the non-hereditary case.

**Proposition 2.** *Let $\mu$ be a Gibbs measure in $\mathcal{G}$. For every bounded non-negative measurable function $f$ from $\mathbb{R}^d \times \mathcal{M}(\mathbb{R}^d)$ to $\mathbb{R}$, we have*

$$\int_{\mathbb{R}^d \times \mathcal{M}(\mathbb{R}^d)} \mathbb{1}_{\mathcal{M}_\infty(\mathbb{R}^d)}(\gamma) f(x, \gamma) \mathcal{C}^!_\mu(\mathrm{d}x, \mathrm{d}\gamma) \qquad (13)$$
$$= \int_{\mathbb{R}^d} \int_{\mathcal{M}(\mathbb{R}^d)} f(x, \gamma) \mathrm{e}^{-h(x,\gamma)} \lambda(\mathrm{d}x) \mu(\mathrm{d}\gamma).$$

We have seen in Remark 2 above that, in general, the support of $\mathcal{C}^!_\mu$ is not included in $\mathbb{R}^d \times \mathcal{M}_\infty(\mathbb{R}^d)$. The Campbell equilibrium equation (13) shows that $\mathcal{C}^!_\mu$ restricted to $\mathbb{R}^d \times \mathcal{M}_\infty(\mathbb{R}^d)$ is absolutely continuous with respect to $\lambda \otimes \mu$ with density $\mathrm{e}^{-h}$.

Let us point out that this proposition is also valid in the non-stationary case. The intensity $\lambda$ is then replaced by any locally finite measure.

Moreover, it is important to note that the converse of Proposition 2 is not true, which means that (13) does not characterize the measure $\mu$. Consider, for example, a measure $\mu$ such that, almost surely, $\gamma$ does not contain any removable points. (13) then becomes the trivial equation $0 = 0$. In fact, the equilibrium equation (13) is interesting only if, $\mu$ almost surely, $\gamma$ contains some removable points.

**Proof of Proposition 2.** Let $\mu$ be a Gibbs measure and $f$ a bounded non-negative measurable function from $\mathbb{R}^d \times \mathcal{M}(\mathbb{R}^d)$ to $\mathbb{R}$. Let $\Lambda$ be a bounded set in $\mathbb{R}^d$. By the definition (5) of the reduced Campbell measure, we have

$$\int_{\mathbb{R}^d \times \mathcal{M}(\mathbb{R}^d)} \mathbb{1}_{\mathcal{M}_\infty(\mathbb{R}^d)}(\gamma) \mathbb{1}_\Lambda(x) f(x,\gamma) \mathcal{C}^!_\mu(\mathrm{d}x, \mathrm{d}\gamma)$$
$$= \int_{\mathcal{M}(\mathbb{R}^d)} \int_{\mathbb{R}^d} \mathbb{1}_{\mathcal{M}_\infty(\mathbb{R}^d)}(\gamma - \delta_x) \mathbb{1}_\Lambda(x) f(x, \gamma - \delta_x) \gamma(\mathrm{d}x) \mu(\mathrm{d}\gamma).$$

By the definition of the local energy, Proposition 1, equation (9) and the DLR equations (11), we have

$$\int_{\mathbb{R}^d \times \mathcal{M}(\mathbb{R}^d)} \mathbb{1}_{\mathcal{M}_\infty(\mathbb{R}^d)}(\gamma) \mathbb{1}_\Lambda(x) f(x,\gamma) \mathcal{C}^!_\mu(\mathrm{d}x, \mathrm{d}\gamma)$$
$$= \int_{\mathcal{M}(\mathbb{R}^d)} \sum_{x \in \mathcal{R}(\gamma) \cap \Lambda} \mathbb{1}_\Lambda(x) f(x, \gamma - \delta_x) \mu(\mathrm{d}\gamma)$$
$$= \int_{\mathcal{M}(\mathbb{R}^d)} \int_{\mathcal{M}(\mathbb{R}^d)} \sum_{x \in \mathcal{R}(\gamma') \cap \Lambda} f(x, \gamma' - \delta_x) \Xi_\Lambda(\mathrm{d}\gamma', \gamma) \mu(\mathrm{d}\gamma)$$
$$= \int_{\mathcal{M}(\mathbb{R}^d)} \int_{\mathcal{M}(\Lambda)} \sum_{x \in \mathcal{R}(\gamma'_\Lambda + \gamma_{\Lambda^c}) \cap \Lambda} f(x, \gamma'_\Lambda + \gamma_{\Lambda^c} - \delta_x) \frac{e^{-H_\Lambda(\gamma'_\Lambda + \gamma_{\Lambda^c})}}{Z_\Lambda(\gamma_{\Lambda^c})} \pi_\Lambda(\mathrm{d}\gamma'_\Lambda) \mu(\mathrm{d}\gamma)$$
$$= \int_{\mathcal{M}(\mathbb{R}^d)} \int_{\mathcal{M}(\Lambda)} \sum_{x \in \mathcal{R}(\gamma'_\Lambda + \gamma_{\Lambda^c}) \cap \Lambda} f(x, \gamma'_\Lambda + \gamma_{\Lambda^c} - \delta_x) e^{-h(x, \gamma'_\Lambda + \gamma_{\Lambda^c} - \delta_x)}$$
$$\times \frac{e^{-H_\Lambda(\gamma'_\Lambda + \gamma_{\Lambda^c} - \delta_x)}}{Z_\Lambda(\gamma_{\Lambda^c})} \pi_\Lambda(\mathrm{d}\gamma'_\Lambda) \mu(\mathrm{d}\gamma)$$
$$= \int_{\mathcal{M}(\mathbb{R}^d)} \int_{\mathcal{M}(\Lambda)} \sum_{x \in \Lambda} \mathbb{1}_{\mathcal{R}(\gamma'_\Lambda + \gamma_{\Lambda^c} + \delta_x - \delta_x)}(x) f(x, \gamma'_\Lambda + \gamma_{\Lambda^c} - \delta_x) e^{-h(x, \gamma'_\Lambda + \gamma_{\Lambda^c} - \delta_x)}$$
$$\times \frac{e^{-H_\Lambda(\gamma'_\Lambda + \gamma_{\Lambda^c} - \delta_x)}}{Z_\Lambda(\gamma_{\Lambda^c})} \pi_\Lambda(\mathrm{d}\gamma'_\Lambda) \mu(\mathrm{d}\gamma)$$
$$= \int_{\mathcal{M}(\mathbb{R}^d)} \int_{\Lambda \times \mathcal{M}(\Lambda)} \mathbb{1}_{\mathcal{R}(\gamma'_\Lambda + \gamma_{\Lambda^c} + \delta_x)}(x) f(x, \gamma'_\Lambda + \gamma_{\Lambda^c}) e^{-h(x, \gamma'_\Lambda + \gamma_{\Lambda^c})}$$
$$\times \frac{e^{-H_\Lambda(\gamma'_\Lambda + \gamma_{\Lambda^c})}}{Z_\Lambda(\gamma_{\Lambda^c})} \mathcal{C}^!_{\pi_\Lambda}(\mathrm{d}x, \mathrm{d}\gamma'_\Lambda) \mu(\mathrm{d}\gamma).$$

The well-known Campbell equilibrium equation for the unit rate Poisson process (that is, $\mathcal{C}^!_{\pi_\Lambda} = \lambda \otimes \pi_\Lambda$) gives

$$\int_{\mathbb{R}^d \times \mathcal{M}(\mathbb{R}^d)} \mathbb{1}_{\mathcal{M}_\infty(\mathbb{R}^d)}(\gamma) \mathbb{1}_\Lambda(x) f(x,\gamma) \mathcal{C}^!_\mu(\mathrm{d}x, \mathrm{d}\gamma)$$

$$= \int_{\mathcal{M}(\mathbb{R}^d)} \int_{\Lambda \times \mathcal{M}(\Lambda)} \mathbb{1}_{\mathcal{R}(\gamma'_\Lambda + \gamma_{\Lambda^c} + \delta_x)}(x) f(x, \gamma'_\Lambda + \gamma_{\Lambda^c}) e^{-h(x, \gamma'_\Lambda + \gamma_{\Lambda^c})}$$

$$\times \frac{e^{-H_\Lambda(\gamma'_\Lambda + \gamma_{\Lambda^c})}}{Z_\Lambda(\gamma_{\Lambda^c})} \lambda \otimes \pi_\Lambda(\mathrm{d}x, \mathrm{d}\gamma'_\Lambda) \mu(\mathrm{d}\gamma).$$

Again using (9), we have

$$\int_{\mathbb{R}^d \times \mathcal{M}(\mathbb{R}^d)} \mathbb{1}_{\mathcal{M}_\infty(\mathbb{R}^d)}(\gamma) \mathbb{1}_\Lambda(x) f(x,\gamma) \mathcal{C}^!_\mu(\mathrm{d}x, \mathrm{d}\gamma)$$

$$= \int_{\mathbb{R}^d} \int_{\mathcal{M}(\mathbb{R}^d)} \int_{\mathcal{M}(\mathbb{R}^d)} \mathbb{1}_{\mathcal{R}(\gamma' + \delta_x) \cap \Lambda}(x) f(x, \gamma') e^{-h(x,\gamma')} \lambda(\mathrm{d}x) \Xi_\Lambda(\mathrm{d}\gamma', \gamma) \mu(\mathrm{d}\gamma).$$

Note that if $\gamma'$ is in $\mathcal{M}_\infty(\mathbb{R}^d)$ and $x$ is in $\mathbb{R}^d$, then $x$ is in $\mathcal{R}(\gamma' + \delta_x)$. Therefore, we have

$$\int_{\mathbb{R}^d \times \mathcal{M}(\mathbb{R}^d)} \mathbb{1}_{\mathcal{M}_\infty(\mathbb{R}^d)}(\gamma) \mathbb{1}_\Lambda(x) f(x,\gamma) \mathcal{C}^!_\mu(\mathrm{d}x, \mathrm{d}\gamma)$$

$$= \int_{\mathbb{R}^d} \int_{\mathcal{M}(\mathbb{R}^d)} \int_{\mathcal{M}(\mathbb{R}^d)} \mathbb{1}_\Lambda(x) f(x, \gamma') e^{-h(x,\gamma')} \lambda(\mathrm{d}x) \Xi_\Lambda(\mathrm{d}\gamma', \gamma) \mu(\mathrm{d}\gamma)$$

$$= \int_{\mathbb{R}^d} \int_{\mathcal{M}(\mathbb{R}^d)} \mathbb{1}_\Lambda(x) f(x, \gamma) e^{-h(x,\gamma)} \lambda(\mathrm{d}x) \mu(\mathrm{d}\gamma).$$

Relation (13) is proved for every $\Lambda$ in $\mathcal{B}(\mathbb{R}^d)$, so it is proved for $\Lambda = \mathbb{R}^d$ as well. □

## 4. Consistency of the pseudo-likelihood estimator

In this section, we deal with the parametric estimation of (non-hereditary) stationary Gibbs measures. Our aim is to prove the asymptotic consistency of the estimation procedure when the observation window of the Gibbs point process increases to $\mathbb{R}^d$. We focus on the pseudo-likelihood procedure, which appears to be asymptotically validated for a large class of interactions. This procedure, quicker in practice than the maximum likelihood approach, may constitute a first consistent estimation in applications and may therefore be used to quickly fit several alternative models.

We suppose that the family of energies $(H_\Lambda)$ depends on a positive parameter $\alpha^*$ and on a multiple parameter $\theta^* = (\theta^*_1, \ldots, \theta^*_p)$. The first parameterizes the support of the energy (that is, when the energy is equal to positive infinity), while the second parameterizes the energy when it is finite. These two parameters play very different roles, as the set of assumptions made in this section will make clear. Note that it would be easy to consider a vectorial hardcore parameter, but, for reasons of clarity, we chose to focus solely on a real hardcore parameter.

Let $(\alpha, \theta)$ be two parameters in $\mathbb{R}^+ \times \Theta$, where $\Theta$ is a bounded open set in $\mathbb{R}^p$. Let us denote by $(H_\Lambda^{\alpha,\theta})_{\Lambda \in \mathcal{B}(\mathbb{R}^d)}$ the parametric family of energies and $\mathcal{G}^{\alpha,\theta}$ the set of stationary Gibbs measures for this family of energies. The following assumption S1 implies that the support of the energy is parameterized by $\alpha$ only, not by $\theta$.

S1. For all $\gamma \in \mathcal{M}(\mathbb{R}^d)$, $\alpha \in \mathbb{R}^+$ and $\theta$, $\theta'$ in $\Theta$,

$$\forall \Lambda \in \mathcal{B}(\mathbb{R}^d) \qquad H_\Lambda^{\alpha,\theta}(\gamma) < \infty \quad \iff \quad H_\Lambda^{\alpha,\theta'}(\gamma) < \infty.$$

This assumption claims that the set of configurations $\gamma$ in $\mathcal{M}(\mathbb{R}^d)$ which have a locally finite energy for the family $(H_\Lambda^{\alpha,\theta})$ depends only on $\alpha$. We therefore denote this set by $\mathcal{M}_\infty^\alpha(\mathbb{R}^d)$. The same remark is true for the set of removable points in $\gamma$, and we similarly denote this set by $\mathcal{R}^\alpha(\gamma)$. Finally, for every $x$ in $\mathcal{R}^\alpha(\gamma)$, we define $h^{\alpha,\theta}(x, \gamma - \delta_x)$, the energy of $x$ in $\gamma - \delta_x$, as in (7).

For all $\alpha \in \mathbb{R}^+$ and $\theta \in \Theta$, we define the pseudo-likelihood function at $\gamma \in \mathcal{M}_\infty^\alpha(\mathbb{R}^d)$ as

$$PLL_{\Lambda_n}(\gamma, \alpha, \theta) = \frac{1}{|\Lambda_n|} \left[ \int_{\Lambda_n} \exp(-h^{\alpha,\theta}(x, \gamma)) \, dx + \sum_{x \in \mathcal{R}^\alpha(\gamma) \cap \Lambda_n} h^{\alpha,\theta}(x, \gamma - \delta_x) \right], \quad (14)$$

where $\Lambda_n$ denotes the domain of observation of the sample and $|\Lambda_n|$ is its Lebesgue measure.

The purpose of this section is twofold. First, when the hardcore parameter is supposed to be known, we aim to prove the consistency of the pseudo-likelihood estimation of $\theta^*$, without assuming the heredity of the interaction. This problem is addressed in Section 4.1. Since we are in a non-hereditary setting, the definition (14) of the pseudo-likelihood contrast function differs from the common one: it involves the set of removable points. In the hereditary case, the proof of the consistency relies mainly on the ergodic theorem and the Nguyen–Zessin equation. In the non-hereditary case, since the Campbell equilibrium equation is modified into (13), we must check carefully the consistency of the pseudo-likelihood procedure. We prove in Theorem 1 that the non-hereditary setting does not actually modify the range of validity of the consistency since the set of hypotheses on the smooth interaction parameter $\theta$ is the same as in the hereditary case (for example, as in [2]).

The second purpose, addressed in Section 4.2, is to prove the consistency of the two-step estimation procedure, without assuming the heredity of the interaction. The support parameter $\alpha^*$ is first estimated and plugged into the pseudo-likelihood contrast function. $\theta^*$ is then estimated by the common pseudo-likelihood procedure. To the best of our knowledge, this two-step procedure, although widely used in practice, has not yet been proven to be consistent, even in the hereditary setting. Apart from the usual conditions on the smooth interaction parameter, some additional assumptions on the support parameter $\alpha$ are needed (see S2–S4): roughly speaking, they suppose that the sets $\mathcal{M}_\infty^\alpha(\mathbb{R}^d)$ are embedded and continuous with respect to $\alpha$.

In Section 5, some examples of models involving both a hardcore and a smooth interaction are presented. We prove that they satisfy all the assumptions of this section.

## 4.1. Consistency of $\hat{\theta}_n$ when the support parameter $\alpha^*$ is known

Let us suppose that the first parameter $\alpha^*$ is known. We estimate $\theta^*$ by

$$\hat{\theta}_n = \underset{\theta \in \Theta}{\operatorname{argmin}}\, PLL_{\Lambda_n}(\gamma, \alpha^*, \theta). \tag{15}$$

To prove the consistency of this estimator, we need the following set of hypotheses:

H1. $(\Lambda_n)_{n\geq 1}$ is an increasing sequence of convex and compact sets such that $|\Lambda_n| \to \mathbb{R}^d$;

H2. the energy function $h$ is invariant by translation, that is, for all $x$ and $y$ in $\mathbb{R}^d$ and all $\gamma \in \mathcal{M}(\mathbb{R}^d)$, $h(x+y, \gamma_y) = h(x,\gamma)$, where $\gamma_y$ is the configuration $\gamma$ translated by $y$;

H3. $\exists K \geq 0$ such that for all $(\alpha,\theta)$, $\gamma \in \mathcal{M}_\infty^\alpha(\mathbb{R}^d)$ and $x \in \mathbb{R}^d$,

$$h^{\alpha,\theta}(x,\gamma) \geq -K;$$

H4. for all $\theta$ in $\Theta$,

$$h^{\alpha^*,\theta}(0,\cdot)\exp(-h^{\alpha^*,\theta^*}(0,\cdot)) \in L^1(\mu^{\alpha^*,\theta^*}),$$

with the convention $\infty \mathrm{e}^{-\infty} = 0$;

H5. for all $\theta$ in $\Theta \setminus \{\theta^*\}$,

$$\mu^{\alpha^*,\theta^*}(h^{\alpha^*,\theta^*}(0,\cdot) \neq h^{\alpha^*,\theta}(0,\cdot)) > 0;$$

H6. one can find a real function $\delta$ with $\delta(x) \to 0$ when $x \to 0$ and $g \in L^1(\mu^{\alpha^*,\theta^*})$, such that $\forall(\theta,\theta') \in \Theta^2$, $\forall \gamma \in \mathcal{M}_\infty^{\alpha^*}(\mathbb{R}^d)$, if $h^{\alpha^*,\theta}(0,\gamma) < +\infty$, then

$$|h^{\alpha^*,\theta}(0,\gamma) - h^{\alpha^*,\theta'}(0,\gamma)| \leq g(\gamma)\delta(|\theta - \theta'|), \qquad \mu^{\alpha^*,\theta^*}\text{-a.e.}$$

H1 is a natural assumption concerning the domain of observation. H2 and H3 state that $h$ is invariant under translation and is locally stable. Local stability is a stronger assumption than the classical stability hypothesis encountered in statistical mechanics, yet it is fulfilled in many Gibbs models (see, for instance, [1, 5]). H4 is a technical integrability assumption which holds in most models (see [1, 5, 14]). H5 guarantees that $\theta$ is a proper parameter for the energy. From H5, we deduce that, $\mu^{\alpha^*,\theta^*}$-a.s., there exists some configuration $\gamma$ such that $h^{\alpha^*,\theta^*}(0,\gamma) < +\infty$. One could then prove, using the ergodic theorem, that $\gamma$ contains almost surely some removable points. As a consequence, the sum involved in (14) is non-empty. Finally, H6 specifies the sense in which $\theta \mapsto h^{\alpha^*,\theta}(0,\gamma)$ is continuous at $\theta^*$.

Most of these hypotheses are similar to the assumptions found in [2]. Indeed, the latter deal with the pseudo-likelihood estimation in the general hereditary case. We adapt their scheme to the non-hereditary setting. The assumptions and the proofs are, therefore, in the same spirit.

**Theorem 1.** *Let $\mu^{\alpha^*,\theta^*} \in \mathcal{G}^{\alpha^*,\theta^*}$. Under* S1 *and* H1–H6*, the estimator $\hat{\theta}_n$ defined by (15) is strongly consistent, that is, $\mu^{\alpha^*,\theta^*}$-a.e.,*

$$\lim_{n\to\infty} \hat{\theta}_n = \theta^*. \tag{16}$$

The pseudo-likelihood procedure is a minimum contrast estimation. This point of view has been used in Jensen and Künsch [9] and in Billiot, Coeurjolly and Drouilhet [2] to prove the consistency and asymptotic normality of their estimator. It relies on Theorem 3.4.3 of Guyon, established in [8].

First, note that we only have to prove Theorem 1 for ergodic measures $\mu^{\alpha^*,\theta^*}$. If $\mu^{\alpha^*,\theta^*}$ is not ergodic, it can be represented as a mixture of ergodic stationary Gibbs measures (see [12]). Therefore, from now on, $\mu^{\alpha^*,\theta^*}$ is assumed to be ergodic.

The following Lemmas 1 and 2 allow us to apply Theorem 3.4.3 in [8], which yields (16).

Let

$$K_n(\theta, \theta^*) = PLL_{\Lambda_n}(\gamma, \alpha^*, \theta) - PLL_{\Lambda_n}(\gamma, \alpha^*, \theta^*). \tag{17}$$

We prove below that $K_n$ is a proper contrast function.

**Lemma 1.** *Under* S1 *and* H1–H5*, for all $\theta \in \Theta$ and $\mu^{\alpha^*,\theta^*}$-a.e.,*

$$\lim_{n\to\infty} K_n(\theta, \theta^*) = K(\theta, \theta^*),$$

*where $K(\cdot, \theta^*)$ is a deterministic positive function which has a unique minimum at $\theta^*$.*

**Proof.** The main point of the proof consists of proving that for all $\theta \in \Theta$ and $\mu^{\alpha^*,\theta^*}$-almost every $\gamma$,

$$\begin{aligned}
\lim_{n\to\infty} &PLL_{\Lambda_n}(\gamma, \alpha^*, \theta) \\
&= E_{\alpha^*,\theta^*}[\exp(-h^{\alpha^*,\theta}(0,\gamma)) + h_{\alpha^*,\theta}(0,\gamma)\exp(-h^{\alpha^*,\theta^*}(0,\gamma))],
\end{aligned} \tag{18}$$

where $E_{\alpha^*,\theta^*}$ denotes the expectation under $\mu^{\alpha^*,\theta^*}$.

Thanks to H1–H3, we can apply the ergodic theorem (see [17]):

$$\lim_{n\to\infty} \frac{1}{|\Lambda_n|} \int_{\Lambda_n} \exp(-h^{\alpha^*,\theta}(x,\gamma))\,\mathrm{d}x = E_{\alpha^*,\theta^*}\left[\int_{[0,1]^d} \exp(-h^{\alpha^*,\theta}(x,\gamma))\,\mathrm{d}x\right].$$

The stationarity of $\mu^{\alpha^*,\theta^*}$ yields

$$\lim_{n\to\infty} \frac{1}{|\Lambda_n|} \int_{\Lambda_n} \exp(-h^{\alpha^*,\theta}(x,\gamma))\,\mathrm{d}x = E_{\alpha^*,\theta^*}[\exp(-h^{\alpha^*,\theta}(0,\gamma))],$$

which proves the first part of (18).

For the second part, note first that from Proposition 2 and the stationarity of $\mu^{\alpha^*,\theta^*}$,

$$\mathcal{C}^!_{\mu^{\alpha^*,\theta^*}}(\mathbb{1}_{\mathcal{M}^{\alpha^*}_\infty(\mathbb{R}^d)}(\gamma)|h_{\alpha^*,\theta}(x,\gamma)|\mathbb{1}_{x\in[0,1]^d}) = E_{\alpha^*,\theta^*}[|h_{\alpha^*,\theta}(0,\gamma)|\mathrm{e}^{-h_{\alpha^*,\theta^*}(0,\gamma)}].$$

From Definition 5 of the reduced Campbell measure, we deduce that

$$E_{\alpha^*,\theta^*}\left[\left|\sum_{x\in\mathcal{R}^{\alpha^*,\theta}(\gamma)\cap[0,1]^d} h^{\alpha^*,\theta}(x,\gamma-\delta_x)\right|\right] \leq E_{\alpha^*,\theta^*}\left[\sum_{x\in\mathcal{R}^{\alpha^*,\theta}(\gamma)\cap[0,1]^d} |h^{\alpha^*,\theta}(x,\gamma-\delta_x)|\right]$$

$$\leq \mathcal{C}^!_{\mu^{\alpha^*,\theta^*}}(\mathbb{1}_{\mathcal{M}^{\alpha^*}_\infty(\mathbb{R}^d)}(\gamma)|h_{\alpha^*,\theta}(x,\gamma)|\mathbb{1}_{x\in[0,1]^d}),$$

which is finite, thanks to H4. Therefore, we can again use the ergodic theorem and $\mu^{\alpha^*,\theta^*}$-a.e.

$$\lim_{n\to\infty}\frac{1}{|\Lambda_n|}\sum_{x\in\mathcal{R}^{\alpha^*}(\gamma)\cap\Lambda_n} h^{\alpha^*,\theta}(x,\gamma-\delta_x) = E_{\alpha^*,\theta^*}\left[\sum_{x\in\mathcal{R}^{\alpha^*}(\gamma)\cap[0,1]^d} h^{\alpha^*,\theta}(x,\gamma-\delta_x)\right]$$

$$= \mathcal{C}^!_{\mu^{\alpha^*,\theta^*}}(\mathbb{1}_{\mathcal{M}^{\alpha^*}_\infty(\mathbb{R}^d)}(\gamma)h_{\alpha^*,\theta}(x,\gamma)\mathbb{1}_{x\in[0,1]^d})$$

$$= E_{\alpha^*,\theta^*}[h_{\alpha^*,\theta}(0,\gamma)\mathrm{e}^{-h_{\alpha^*,\theta^*}(0,\gamma)}].$$

The convergence (18) is proved.

From the definition (17) of $K_n$, we deduce from (18) that $\mu^{\alpha^*,\theta^*}$-a.s.,

$$\lim_{n\to\infty} K_n(\theta,\theta^*)$$
$$= E_{\alpha^*,\theta^*}[\mathrm{e}^{-h_{\alpha^*,\theta^*}(0,\gamma)}(\mathrm{e}^{h_{\alpha^*,\theta^*}(0,\gamma)-h_{\alpha^*,\theta}(0,\gamma)} - [h_{\alpha^*,\theta^*}(0,\gamma) - h_{\alpha^*,\theta}(0,\gamma)] - 1)]$$
$$:= K(\theta,\theta^*).$$

According to H5 and the behavior of $t \mapsto \mathrm{e}^t - t - 1$, we conclude that $K$ is positive and vanishes if and only if $\theta = \theta^*$. □

Let us define the modulus of continuity of $PLL_{\Lambda_n}(\gamma,\alpha^*,\cdot)$ by

$$W_n(\eta) = \sup_{|\theta-\theta'|\leq\eta}\{|PLL_{\Lambda_n}(\gamma,\alpha^*,\theta) - PLL_{\Lambda_n}(\gamma,\alpha^*,\theta')|\}.$$

**Lemma 2.** *Under* S1 *and* H1–H6, *there exists a positive sequence* $(\epsilon_k)_{k\geq 1}$ *with* $\epsilon_k \to 0$ *when* $k \to \infty$, *such that for all* $k \geq 1$,

$$\mu^{\alpha^*,\theta^*}\left(\limsup_{n\to\infty}\left(W_n\left(\frac{1}{k}\right)\geq\epsilon_k\right)\right) = 0.$$

*As a consequence, the functions* $\theta \mapsto PLL_{\Lambda_n}(\gamma,\alpha^*,\theta)$ *and* $\theta \mapsto K_n(\theta,\theta^*)$ *are continuous.*

**Proof.**
$$W_n\left(\frac{1}{k}\right) \leq W_{1,n}\left(\frac{1}{k}\right) + W_{2,n}\left(\frac{1}{k}\right),$$

where
$$W_{1,n}\left(\frac{1}{k}\right) = \sup_{|\theta - \theta'| \leq 1/k} \left\{ \left| \frac{1}{|\Lambda_n|} \int_{\Lambda_n} \exp(-h^{\alpha^*,\theta}(x,\gamma)) - \exp(-h^{\alpha^*,\theta'}(x,\gamma)) \, dx \right| \right\}$$

and
$$W_{2,n}\left(\frac{1}{k}\right) = \sup_{|\theta - \theta'| \leq 1/k} \left\{ \left| \frac{1}{|\Lambda_n|} \sum_{x \in \mathcal{R}^{\alpha^*}(\gamma) \cap \Lambda_n} h^{\alpha^*,\theta}(x,\gamma - \delta_x) - h^{\alpha^*,\theta'}(x,\gamma - \delta_x) \right| \right\}.$$

From S1, we know that $h^{\alpha^*,\theta}(x,\gamma)$ and $h^{\alpha^*,\theta'}(x,\gamma)$ are either both finite or both infinite. Thus, from H6 and H3,

$$W_{1,n}\left(\frac{1}{k}\right) \leq \frac{e^K}{|\Lambda_n|} \int_{\Lambda_n} \delta\left(\frac{1}{k}\right) |g(\gamma_{-x})| \, dx,$$

where $\gamma_{-x}$, for $n$ sufficiently large denotes the configuration $\gamma$ translated by $-x$. Now, according to the ergodic theorem and by the stationarity of $\mu^{\alpha^*,\theta^*}$,

$$W_{1,n}\left(\frac{1}{k}\right) \leq 2\delta\left(\frac{1}{k}\right) e^K E_{\alpha^*,\theta^*}[|g(\gamma)|].$$

From Proposition 1 and S1, $h^{\alpha^*,\theta}(x,\gamma - \delta_x)$ and $h^{\alpha^*,\theta'}(x,\gamma - \delta_x)$ are both finite when $x \in \mathbb{R}^{\alpha^*}(\gamma)$. Thus, from H6 and H3,

$$W_{2,n}\left(\frac{1}{k}\right) \leq \frac{1}{|\Lambda_n|} \sum_{x \in \mathcal{R}^{\alpha^*}(\gamma) \cap \Lambda_n} \delta\left(\frac{1}{k}\right) |g((\gamma - \delta_x)_{-x})|.$$

Applying the ergodic theorem to the right-hand side, Proposition 2 and the stationarity of $\mu^{\alpha^*,\theta^*}$ then lead to

$$W_{2,n}\left(\frac{1}{k}\right) \leq 2\delta\left(\frac{1}{k}\right) E_{\alpha^*,\theta^*}[|g(\gamma)|e^{-h^{\alpha^*,\theta^*}(0,\gamma)}] \leq 2e^K \delta\left(\frac{1}{k}\right) E_{\alpha^*,\theta^*}[|g(\gamma)|].$$

Therefore, for $n$ sufficiently large

$$W_n\left(\frac{1}{k}\right) \leq c\delta\left(\frac{1}{k}\right) E_{\alpha^*,\theta^*}(|g(\gamma)|),$$

where $c$ is a positive constant. Finally,

$$\mu^{\alpha^*,\theta^*}\left(\limsup_{n \to \infty} \left(W_n\left(\frac{1}{k}\right) \geq \epsilon_k\right)\right) \leq \mu^{\alpha^*,\theta^*}\left(c\delta\left(\frac{1}{k}\right) E_{\alpha^*,\theta^*}(|g(\gamma)|) \geq \epsilon_k\right),$$

which vanishes if one chooses, for instance, $\epsilon_k = 2c\delta(\frac{1}{k})E_{\alpha^*,\theta^*}(|g(\gamma)|)$. □

## 4.2. Consistency of $(\hat{\alpha}_n, \hat{\theta}_n)$ when the support parameter $\alpha^*$ is unknown

In this section, the two-step estimation procedure is considered. The consistency of this procedure requires some additional assumptions concerning the support parameter $\alpha$. It must be remembered that $\alpha$ is a positive parameter satisfying the support hypothesis S1. Moreover, we assume S2–S4 below.

S2. For all $\gamma \in \mathcal{M}(\mathbb{R}^d)$, $\alpha$ and $\alpha'$ in $\mathbb{R}^+$, $\theta \in \Theta$ and $\Lambda \in \mathcal{B}(\mathbb{R}^d)$,

$$\alpha \leq \alpha' \implies [H_\Lambda^{\alpha,\theta}(\gamma) < \infty \Rightarrow H_\Lambda^{\alpha',\theta}(\gamma) < \infty].$$

S3. Letting $\alpha < \alpha^*$, there exists $R_\alpha > 0$ such that for every Gibbs measure $\mu^{\alpha^*,\theta^*}$ in $\mathcal{G}^{\alpha^*,\theta^*}$,

$$\mu^{\alpha^*,\theta^*}(H_{\mathcal{B}(0,R_\alpha)}^{\alpha,\theta^*}(\gamma) = +\infty) > 0.$$

S4. For all $\gamma \in \mathcal{M}(\mathbb{R}^d)$, $\Lambda \in \mathcal{B}(\mathbb{R}^d)$ in $\mathbb{R}^+$, $\theta$ in $\Theta$ and $\alpha' > 0$ such that $H_\Lambda^{\alpha',\theta}(\gamma) < +\infty$, there exists $\alpha < \alpha'$ such that $H_\Lambda^{\alpha,\theta}(\gamma) < +\infty$.

Assumption S2 implies that if $\alpha \leq \alpha'$, then $\mathcal{M}_\infty^\alpha(\mathbb{R}^d) \subset \mathcal{M}_\infty^{\alpha'}(\mathbb{R}^d)$ and $\mathcal{R}^\alpha(\gamma) \subset \mathcal{R}^{\alpha'}(\gamma)$. As a consequence, assumption S3 appears to be a natural assumption concerning the support parameter. Roughly speaking, it claims that for an underestimated support of the energies, one could encounter some forbidden configurations with a non-negligible probability.

For a given $\gamma \in \mathcal{M}(\mathbb{R}^d)$, S4 may be viewed as an assumption of continuity of $\mathcal{R}^\alpha(\gamma)$ with respect to $\alpha$. Indeed, it implies that

$$\mathbb{R}^{\alpha'}(\gamma) = \bigcup_{\alpha < \alpha'} \mathbb{R}^\alpha(\gamma). \tag{19}$$

Nevertheless, in general, $\mathcal{M}_\infty^{\alpha'}(\mathbb{R}^d)$ is not equal to $\bigcup_{\alpha<\alpha'} \mathcal{M}_\infty^\alpha(\mathbb{R}^d)$. Some consequences of S4 are proved in the following lemma.

**Lemma 3.** *Under the assumption* S2 *and* S4*, we have*

$$\lim_{\epsilon \to 0^+} \mu^{\alpha^*,\theta^*}(0 \notin \mathbb{R}^{\alpha^*-\epsilon}(\gamma + \delta_0)) = 0, \tag{20}$$

$$\lim_{\epsilon \to 0^+} \mu^{\alpha^*,\theta^*}(h^{\alpha^*-\epsilon,\theta}(0,\gamma) = +\infty \text{ and} \\ h^{\alpha^*,\theta}(0,\gamma) < \infty \text{ and } 0 \in \mathbb{R}^{\alpha^*-\epsilon}(\gamma + \delta_0)) = 0. \tag{21}$$

The proof of this lemma is relegated to Section 4.2.3.

*4.2.1. Consistency of $\hat{\alpha}_n$*

According to S2, we estimate the support parameter $\alpha^*$ by the natural estimator

$$\hat{\alpha}_n = \inf\{\alpha > 0, H_{\Lambda_n}^{\alpha,\theta}(\gamma) < \infty\}. \tag{22}$$

Let us note that, from S1, $\hat{\alpha}_n$ does not depend on $\theta$. Moreover, it is well defined since $H_{\Lambda_n}^{\alpha^*,\theta}(\gamma) < \infty$.

**Proposition 3.** *Under* S1–S3 *and* H1–H3, $\mu^{\alpha^*,\theta^*}$-*a.s.,*

$$\lim_{n \to \infty} \hat{\alpha}_n = \alpha^*.$$

**Proof.** From H1 and relation (4) concerning the energies, if $n < m$, then

$$\{\alpha, H_{\Lambda_m}^{\alpha,\theta}(\gamma) < \infty\} \subset \{\alpha, H_{\Lambda_n}^{\alpha,\theta}(\gamma) < \infty\}.$$

Hence, $(\hat{\alpha}_n)$ is an increasing sequence. From (22), it is clear that $\hat{\alpha}_n \leq \alpha^*$. Therefore, $\hat{\alpha}_n \to \tilde{\alpha}$, where $\tilde{\alpha} = \sup_n \hat{\alpha}_n \leq \alpha^*$. Let us prove that $\tilde{\alpha} < \alpha^*$ is impossible.

Let us assume $\tilde{\alpha} < \alpha^*$ and let $\tilde{\alpha}'$ be such that $\tilde{\alpha} < \tilde{\alpha}' < \alpha^*$. Consider the average

$$\frac{1}{|\Lambda_n|} \int_{\Lambda_n} \mathbb{1}_{\{H_{\mathcal{B}(x,R_{\tilde{\alpha}'})}^{\tilde{\alpha}',\theta}(\gamma) = +\infty\}} \, \mathrm{d}x,$$

where $R_{\tilde{\alpha}'}$ is defined in S3. The ergodic theorem applies and for $\mu^{\alpha^*,\theta^*}$-almost every $\gamma$,

$$\lim_{n \to \infty} \frac{1}{|\Lambda_n|} \int_{\Lambda_n} \mathbb{1}_{\{H_{\mathcal{B}(x,R_{\tilde{\alpha}'})}^{\tilde{\alpha}',\theta}(\gamma) = +\infty\}} \, \mathrm{d}x = \mu^{\alpha^*,\theta^*}(H_{\mathcal{B}(0,R_{\tilde{\alpha}'})}^{\tilde{\alpha}',\theta}(\gamma) = +\infty).$$

From S3, this last term is positive, hence for $\mu^{\alpha^*,\theta^*}$-almost every $\gamma$ and for $n$ large enough, there exists $x \in \Lambda_n$ such that the energy $H_{\mathcal{B}(x,R_{\tilde{\alpha}'})}^{\tilde{\alpha}',\theta}(\gamma)$ is not finite. From relation (4), this means that for a sufficiently large $n_0$, $H_{\Lambda_{n_0}}^{\tilde{\alpha}',\theta}(\gamma) = +\infty$. From (22), this implies that $\tilde{\alpha}' \leq \hat{\alpha}_{n_0}$. But $\tilde{\alpha}' > \tilde{\alpha} \geq \hat{\alpha}_{n_0}$. We have a contradiction and so, finally, $\mu^{\alpha^*,\theta^*}$-a.s., $\tilde{\alpha} = \alpha^*$. □

*4.2.2. Consistency of $(\hat{\alpha}_n, \hat{\theta}_n)$*

The parameter $\theta^*$ is estimated, as in Section 4.1, via the pseudo-likelihood procedure. However, $\alpha^*$ is not known and we have to plug $\hat{\alpha}_n$ into definition (14) of $PLL_{\Lambda_n}$. The estimator $\hat{\theta}_n$ is thus defined as

$$\hat{\theta}_n = \underset{\theta \in \Theta}{\operatorname{argmin}} \, PLL_{\Lambda_n}(\gamma, \hat{\alpha}_n, \theta), \tag{23}$$

where $\hat{\alpha}_n$ is the estimator (22) and

$$PLL_{\Lambda_n}(\gamma, \hat{\alpha}_n, \theta) = \frac{1}{\Lambda_n} \left[ \int_{\Lambda_n} \exp(-h^{\hat{\alpha}_n,\theta}(x,\gamma)) \, \mathrm{d}x + \sum_{x \in \mathbb{R}^{\hat{\alpha}_n}(\gamma) \cap \Lambda_n} h^{\hat{\alpha}_n,\theta}(x, \gamma - \delta_x) \right].$$

**Remark 3.** In the expression above, $h^{\hat\alpha_n,\theta}(x,\gamma)$ is well defined only if $H^{\hat\alpha_n,\theta}_{\Lambda_n}(\gamma)<\infty$. From (22), this is not necessarily the case. However, one can consider a new estimator defined as $\tilde\alpha_n = \hat\alpha_n + \epsilon_n$, where $\epsilon_n$ is any positive sequence which asymptotically vanishes. This new estimator has the same $\mu^{\alpha^*,\theta^*}$-a.e. asymptotic properties as $\hat\alpha_n$ and $h^{\tilde\alpha_n,\theta}(x,\gamma)$ is obviously well defined. In the sequel, we continue to denote by $\hat\alpha_n$ the estimator of $\alpha^*$ so that $h^{\hat\alpha_n,\theta}(x,\gamma)$ is always assumed to be well defined. In addition, if $x \in \mathcal{R}^{\hat\alpha_n}(\gamma)$, $h^{\hat\alpha_n,\theta}(x,\gamma-\delta_x)$ is always well defined and, moreover, it is $\mu^{\alpha^*,\theta^*}$-a.e. finite since $H^{\hat\alpha_n,\theta}_{\Lambda_n}(\gamma)<\infty$.

To prove the consistency of $\hat\theta_n$, we need the hypothesis H1–H5 and the following modification of H6.

**H6′.** One can find a real function $\delta$ with $\delta(x) \to 0$ when $x \to 0$ and $g \in L^1(\mu^{\alpha^*,\theta^*})$, such that for all $\alpha \leq \alpha^*$, $(\theta,\theta') \in \Theta^2$ and $\gamma \in \mathcal{M}^{\alpha^*}_\infty(\mathbb{R}^d)$, if $0 \in \mathcal{R}^\alpha(\gamma+\delta_0)$ and $h^{\alpha,\theta}(0,\gamma)<\infty$, then

$$|h^{\alpha,\theta}(0,\gamma) - h^{\alpha^*,\theta'}(0,\gamma)| \leq g(\gamma)[\delta(|\alpha-\alpha^*|) + \delta(|\theta-\theta'|)], \qquad \mu^{\alpha^*,\theta^*}\text{-a.e.}$$

Let us note that in H6′, we must assume that $0 \in \mathcal{R}^\alpha(\gamma+\delta_0)$, which ensures that $h^{\alpha,\theta}(0,\gamma)$ exists when $\gamma \in \mathcal{M}^{\alpha^*}_\infty(\mathbb{R}^d)$. This is not a restriction in our case since we shall apply H6′ with $\alpha = \hat\alpha_n$ (see Remark 3). Besides, from S1, the condition $h^{\alpha,\theta}(0,\gamma)<\infty$ implies that $h^{\alpha^*,\theta}(0,\gamma)<\infty$ as well.

**Theorem 2.** *Let $\mu^{\alpha^*,\theta^*} \in \mathcal{G}^{\alpha^*,\theta^*}$. Under S1–S4, H1–H5 and H6′, the estimators $\hat\alpha_n$ and $\hat\theta_n$ defined by (22) and (23), respectively, are strongly consistent, that is, $\mu^{\alpha^*,\theta^*}$-a.e.,*

$$\lim_{n\to\infty}(\hat\alpha_n,\hat\theta_n) = (\alpha^*,\theta^*). \qquad (24)$$

The main point of the proof of Theorem 2 is the following lemma, which is proved in Section 4.2.3.

**Lemma 4.** *Under S1–S4, H1–H4 and H6′, for all $\theta \in \Theta$,*

$$\lim_{n\to\infty}(PLL_{\Lambda_n}(\gamma,\hat\alpha_n,\theta) - PLL_{\Lambda_n}(\gamma,\alpha^*,\theta)) = 0, \qquad \mu^{\alpha^*,\theta^*}\text{-a.e.}$$

**Proof of Theorem 2.** The proof follows along the same lines as the proof of Theorem 1. Let

$$K'_n(\theta,\theta^*) = PLL_{\Lambda_n}(\gamma,\hat\alpha_n,\theta) - PLL_{\Lambda_n}(\gamma,\hat\alpha_n,\theta^*). \qquad (25)$$

Thanks to Lemmas 1 and 4, it is clear that

$$\lim_{n\to\infty} K'_n(\theta,\theta^*) = K(\theta,\theta^*),$$

where $K(\cdot,\theta^*)$ is the same function as in Lemma 1. Therefore, $K'_n$ is a proper contrast function.

Now, let the modulus of continuity of $PLL_{\Lambda_n}(\gamma, \hat{\alpha}_n, \cdot)$ be

$$W'_n(\eta) = \sup_{|\theta - \theta'| \leq \eta} \{|PLL_{\Lambda_n}(\gamma, \hat{\alpha}_n, \theta) - PLL_{\Lambda_n}(\gamma, \hat{\alpha}_n, \theta')|\}.$$

The result stated in Lemma 2 for $W_n$ is still true for $W'_n$. To prove this, it is sufficient to plug $\hat{\alpha}_n$, instead of $\alpha^*$, into its demonstration. Hence, there exists $\epsilon_k \to 0$ such that, for all $k \geq 1$,

$$\mu^{\alpha^*,\theta^*}\left(\limsup_{n \to \infty} \left(W'_n\left(\frac{1}{k}\right) \geq \epsilon_k\right)\right) = 0.$$

Finally, Theorem 3.4.3 of [8] can be applied and $\hat{\theta}_n$ defined by (23) converges $\mu^{\alpha^*,\theta^*}$-a.e. to $\theta$. This, together with Proposition 3, completes the proof. □

*4.2.3. Proofs of Lemmas 3 and 4.*

**Proof of Lemma 3.** Let us begin with the limit (20). From (19) and S2, we have

$$\lim_{\epsilon \to 0^+} \mu^{\alpha^*,\theta^*}(0 \notin \mathcal{R}^{\alpha^*-\epsilon}(\gamma + \delta_0))$$

$$= \mu^{\alpha^*,\theta^*}\left(\bigcap_{\varepsilon > 0}\{0 \notin \mathcal{R}^{\alpha^*-\epsilon}(\gamma + \delta_0)\}\right) = \mu^{\alpha^*,\theta^*}\left(\overline{\bigcup_{\varepsilon > 0}\{0 \in \mathcal{R}^{\alpha^*-\epsilon}(\gamma + \delta_0)\}}\right)$$

$$= \mu^{\alpha^*,\theta^*}(0 \notin \mathcal{R}^{\alpha^*}(\gamma + \delta_0)) = 0.$$

Now, for the limit (21),

$$\limsup_{\epsilon \to 0^+} \mu^{\alpha^*,\theta^*}(h^{\alpha^*-\epsilon,\theta}(0,\gamma) = +\infty \text{ and } h^{\alpha^*,\theta}(0,\gamma) < \infty \text{ and } 0 \in \mathcal{R}^{\alpha^*-\epsilon}(\gamma + \delta_0))$$

$$\leq \limsup_{\epsilon \to 0^+} \mu^{\alpha^*,\theta^*}(\exists \Lambda \in \mathcal{B}(\mathbb{R}^d) \text{ s.t. } 0 \in \Lambda, H_\Lambda^{\alpha^*-\epsilon,\theta}(\gamma + \delta_0) = +\infty$$

$$\text{and } H_\Lambda^{\alpha^*,\theta}(\gamma + \delta_0) < \infty \text{ and } H_\Lambda^{\alpha^*-\epsilon,\theta}(\gamma) < +\infty)$$

$$\leq \limsup_{\epsilon \to 0^+} \mu^{\alpha^*,\theta^*}(\exists \Lambda \in \mathcal{B}(\mathbb{R}^d) \text{ s.t. } 0 \in \Lambda, H_{\Lambda \cap B(0,1)}^{\alpha^*-\epsilon,\theta}(\gamma + \delta_0) = +\infty$$

$$\text{and } \gamma + \delta_0 \in \mathcal{M}_\infty^{\alpha^*}(\mathbb{R}^d)).$$

Thanks to a monotone class argument,

$$\limsup_{\epsilon \to 0^+} \mu^{\alpha^*,\theta^*}(h^{\alpha^*-\epsilon,\theta}(0,\gamma) = +\infty \text{ and } h^{\alpha^*,\theta}(0,\gamma) < \infty \text{ and } 0 \in \mathcal{R}^{\alpha^*-\epsilon}(\gamma + \delta_0))$$

$$\leq \mu^{\alpha^*,\theta^*}\left(\bigcap_{\epsilon > 0}\{\exists \Lambda \in \mathcal{B}(\mathbb{R}^d) \text{ s.t. } 0 \in \Lambda, H_{\Lambda \cap B(0,1)}^{\alpha^*-\epsilon,\theta}(\gamma + \delta_0) = +\infty \text{ and } \gamma + \delta_0 \in \mathcal{M}_\infty^{\alpha^*}(\mathbb{R}^d)\}\right)$$

$$\leq \mu^{\alpha^*,\theta^*}\left(\{\gamma+\delta_0\in\mathcal{M}_\infty^{\alpha^*}(\mathbb{R}^d)\}\cap\overline{\bigcup_{\epsilon>0}\bigcap_{\substack{\Lambda\in B(\mathbb{R}^d)\\0\in\Lambda}}\{H_{\Lambda\cap B(0,1)}^{\alpha^*-\epsilon,\theta}(\gamma+\delta_0)<+\infty\}}\right)$$

$$\leq \mu^{\alpha^*,\theta^*}\left(\{\gamma+\delta_0\in\mathcal{M}_\infty^{\alpha^*}(\mathbb{R}^d)\}\cap\overline{\bigcup_{\epsilon>0}\{H_{B(0,1)}^{\alpha^*-\epsilon,\theta}(\gamma+\delta_0)<+\infty\}}\right).$$

Finally, from S4,

$$\limsup_{\epsilon\to 0^+}\mu^{\alpha^*,\theta^*}(h^{\alpha^*-\epsilon,\theta}(0,\gamma)=+\infty \text{ and } h^{\alpha^*,\theta}(0,\gamma)<\infty \text{ and } 0\in\mathcal{R}^{\alpha^*-\epsilon}(\gamma+\delta_0))$$

$$\leq \mu^{\alpha^*,\theta^*}(\gamma+\delta_0\in\mathcal{M}_\infty^{\alpha^*}(\mathbb{R}^d) \text{ and } H_{B(0,1)}^{\alpha^*,\theta}(\gamma+\delta_0)=+\infty)=0.$$

Lemma 3 is thus proved. □

**Proof of Lemmas 4.** Let us split the difference as

$$PLL_{\Lambda_n}(\gamma,\hat{\alpha}_n,\theta)-PLL_{\Lambda_n}(\gamma,\alpha^*,\theta)=D_{1,n}+D_{2,n}, \qquad (26)$$

where

$$D_{1,n}=\frac{1}{|\Lambda_n|}\int_{\Lambda_n}\exp(-h^{\hat{\alpha}_n,\theta}(x,\gamma))-\exp(-h^{\alpha^*,\theta}(x,\gamma))\,\mathrm{d}x$$

and

$$D_{2,n}=\frac{1}{|\Lambda_n|}\sum_{x\in\gamma\cap\Lambda_n}(\mathbb{1}_{\mathcal{R}^{\hat{\alpha}_n}(\gamma)}(x)h^{\hat{\alpha}_n,\theta}(x,\gamma-\delta_x)-\mathbb{1}_{\mathcal{R}^{\alpha^*}(\gamma)}(x)h^{\alpha^*,\theta}(x,\gamma-\delta_x)).$$

In the integral of $D_{1,n}$, for a point $x$ belonging to $\Lambda_n$, there are several exclusive cases. The first is $h^{\hat{\alpha}_n,\theta}(x,\gamma)=+\infty$ and $h^{\alpha^*,\theta}(x,\gamma)=+\infty$. In this case, each term vanishes. The second is $h^{\hat{\alpha}_n,\theta}(x,\gamma)<+\infty$ and $h^{\alpha^*,\theta}(x,\gamma)<+\infty$; let us denote by $\Lambda_{1,n}$ the set of such $x$'s. The last case is $h^{\hat{\alpha}_n,\theta}(x,\gamma)=+\infty$ and $h^{\alpha^*,\theta}(x,\gamma)<+\infty$; let us denote by $\Lambda_{2,n}$ the set of such $x$'s. Note that because of S2, $h^{\hat{\alpha}_n,\theta}(x,\gamma)<+\infty$ and $h^{\alpha^*,\theta}(x,\gamma)=+\infty$ is impossible since $\hat{\alpha}_n\leq\alpha^*$. Thus,

$$|D_{1,n}|\leq\frac{1}{|\Lambda_n|}\int_{\Lambda_{1,n}}|e^{-h^{\hat{\alpha}_n,\theta}(x,\gamma)}-e^{-h^{\alpha^*,\theta}(x,\gamma)}|\,\mathrm{d}x+\frac{1}{|\Lambda_n|}\int_{\Lambda_{2,n}}e^{-h^{\alpha^*,\theta}(x,\gamma)}\,\mathrm{d}x.$$

According to H3 and H6',

$$\frac{1}{|\Lambda_n|}\int_{\Lambda_{1,n}}|e^{-h^{\hat{\alpha}_n,\theta}(x,\gamma)}-e^{-h^{\alpha^*,\theta}(x,\gamma)}|\,\mathrm{d}x\leq e^K\frac{1}{|\Lambda_n|}\int_{\Lambda_{1,n}}|g(\gamma_{-x})|\delta(|\hat{\alpha}_n-\alpha^*|)\,\mathrm{d}x,$$

where $\gamma_{-x}$ denotes the configuration $\gamma$ translated by $-x$. If we let $\epsilon>0$, then for $n$ sufficiently large, from Proposition 3,

$$\frac{1}{|\Lambda_n|}\int_{\Lambda_{1,n}}|e^{-h^{\hat{\alpha}_n,\theta}(x,\gamma)}-e^{-h^{\alpha^*,\theta}(x,\gamma)}|\,\mathrm{d}x\leq\epsilon\frac{e^K}{|\Lambda_n|}\int_{\Lambda_n}|g(\gamma_{-x})|\,\mathrm{d}x.$$

Since $g \in L^1(\mu^{\alpha^*,\theta^*})$, the ergodic theorem applies to the average in the right-hand side and, for $n$ sufficiently large,

$$\frac{1}{|\Lambda_n|} \int_{\Lambda_{1,n}} |\mathrm{e}^{-h^{\hat{\alpha}_n,\theta}(x,\gamma)} - \mathrm{e}^{-h^{\alpha^*,\theta}(x,\gamma)}| \,\mathrm{d}x \leq 2\epsilon \mathrm{e}^K E_{\alpha^*,\theta^*}\left(\int_{[0,1]^d} |g(\gamma_{-x})| \,\mathrm{d}x\right).$$

The stationarity of $\mu^{\alpha^*,\theta^*}$ leads to

$$\frac{1}{|\Lambda_n|} \int_{\Lambda_{1,n}} |\mathrm{e}^{-h^{\hat{\alpha}_n,\theta}(x,\gamma)} - \mathrm{e}^{-h^{\alpha^*,\theta}(x,\gamma)}| \,\mathrm{d}x \leq 2\epsilon \mathrm{e}^K E_{\alpha^*,\theta^*}(|g(\gamma)|). \qquad (27)$$

Besides, from H3 and the definition of $\Lambda_{2,n}$,

$$\frac{1}{|\Lambda_n|} \int_{\Lambda_{2,n}} \mathrm{e}^{-h^{\alpha^*,\theta}(x,\gamma)} \,\mathrm{d}x \leq \frac{\mathrm{e}^K}{|\Lambda_n|} \int_{\Lambda_n} \mathbb{1}_{\{h^{\alpha^*,\theta}(x,\gamma)<\infty\}} \mathbb{1}_{\{h^{\hat{\alpha}_n,\theta}(x,\gamma)=+\infty\}} \,\mathrm{d}x.$$

If we let $\epsilon > 0$, then for $n$ sufficiently large, $\hat{\alpha}_n > \alpha^* - \epsilon$. Hence, provided $x \in \mathcal{R}^{\alpha^*-\epsilon}(\gamma + \delta_x)$, we deduce from S1 that $h^{\hat{\alpha}_n,\theta}(x,\gamma) = +\infty$ yields $h^{\alpha^*-\epsilon,\theta}(x,\gamma) = +\infty$. Therefore, for $n$ sufficiently large,

$$\frac{1}{|\Lambda_n|} \int_{\Lambda_{2,n}} \mathrm{e}^{-h^{\alpha^*,\theta}(x,\gamma)} \,\mathrm{d}x$$
$$\leq \frac{\mathrm{e}^K}{|\Lambda_n|} \int_{\Lambda_n} \mathbb{1}_{\{h^{\alpha^*,\theta}(x,\gamma)<\infty\} \cap \{x \notin \mathcal{R}^{\alpha^*-\epsilon}(\gamma+\delta_x)\}} \,\mathrm{d}x$$
$$+ \frac{\mathrm{e}^K}{|\Lambda_n|} \int_{\Lambda_n} \mathbb{1}_{\{h^{\alpha^*,\theta}(x,\gamma)<\infty\} \cap \{x \in \mathcal{R}^{\alpha^*-\epsilon}(\gamma+\delta_x)\} \cap \{h^{\alpha^*-\epsilon,\theta}(x,\gamma)=+\infty\}} \,\mathrm{d}x.$$

According to the ergodic theorem and by the stationarity of $\mu^{\alpha^*,\theta^*}$, for $n$ sufficiently large,

$$\frac{1}{|\Lambda_n|} \int_{\Lambda_{2,n}} \mathrm{e}^{-h^{\alpha^*,\theta}(x,\gamma)} \,\mathrm{d}x \leq 2\mathrm{e}^K \mu^{\alpha^*,\theta^*}(\{h^{\alpha^*,\theta}(0,\gamma) < \infty\} \cap \{0 \notin \mathcal{R}^{\alpha^*-\epsilon}(\gamma+\delta_0)\})$$
$$+ 2\mathrm{e}^K \mu^{\alpha^*,\theta^*}(\{h^{\alpha^*,\theta}(0,\gamma) < \infty\} \cap \{0 \in \mathcal{R}^{\alpha^*-\epsilon}(\gamma+\delta_0)\}$$
$$\cap \{h^{\alpha^*-\epsilon,\theta}(0,\gamma) = +\infty\}),$$

which is less (up to $2\mathrm{e}^K$) than

$$\mu^{\alpha^*,\theta^*}(0 \notin \mathcal{R}^{\alpha^*-\epsilon}(\gamma+\delta_0))$$
$$+ \mu^{\alpha^*,\theta^*}(h^{\alpha^*-\epsilon,\theta}(0,\gamma) = +\infty | \{h^{\alpha^*,\theta}(0,\gamma) < \infty\} \cap \{0 \in \mathcal{R}^{\alpha^*-\epsilon}(\gamma+\delta_0)\}).$$

This last term vanishes when $\epsilon \to 0$, as stated in Lemma 3. This result, together with (27), proves that in (26), $D_{1,n}$ vanishes when $n$ goes to $+\infty$.

Now let us investigate the behavior of $D_{2,n}$ in (26). We have

$$|D_{2,n}| \leq D_{21,n} + D_{22,n},$$

where

$$D_{21,n} = \frac{1}{|\Lambda_n|} \sum_{x \in \gamma \cap \Lambda_n} (\mathbb{1}_{\mathcal{R}^{\alpha^*}(\gamma)}(x) - \mathbb{1}_{\mathcal{R}^{\hat{\alpha}_n}(\gamma)}(x)) |h^{\alpha^*,\theta}(x, \gamma - \delta_x)|$$

and

$$D_{22,n} = \frac{1}{|\Lambda_n|} \sum_{x \in \gamma \cap \Lambda_n} \mathbb{1}_{\mathcal{R}^{\hat{\alpha}_n}(\gamma)}(x) |h^{\hat{\alpha}_n,\theta}(x, \gamma - \delta_x) - h^{\alpha^*,\theta}(x, \gamma - \delta_x)|.$$

If we let $\epsilon > 0$, then for $n$ sufficiently large, $\hat{\alpha}_n > \alpha^* - \epsilon$ and, according to S2, $\mathcal{R}^{\alpha^* - \epsilon}(\gamma) \subset \mathcal{R}^{\hat{\alpha}_n}(\gamma)$. Thus,

$$D_{21,n} \leq \frac{1}{|\Lambda_n|} \sum_{x \in \gamma \cap \Lambda_n} (\mathbb{1}_{\mathcal{R}^{\alpha^*}(\gamma)}(x) - \mathbb{1}_{\mathcal{R}^{\alpha^* - \epsilon}(\gamma)}(x)) |h^{\alpha^*,\theta}(x, \gamma - \delta_x)|.$$

The application to the right-hand side of the ergodic theorem, combined with Proposition 2 and the stationarity of $\mu^{\alpha^*,\theta^*}$, leads to

$$D_{21,n} \leq 2 E_{\alpha^*,\theta^*}[(1 - \mathbb{1}_{\mathcal{R}^{\alpha^* - \epsilon}(\gamma + \delta_0)}(0)) |h^{\alpha^*,\theta}(0,\gamma)| e^{-h^{\alpha^*,\theta^*}(0,\gamma)}]$$

$$\leq 2 E_{\alpha^*,\theta^*}[\mathbb{1}_{0 \notin \mathcal{R}^{\alpha^* - \epsilon}(\gamma + \delta_0)} |h^{\alpha^*,\theta}(0,\gamma)| e^{-h^{\alpha^*,\theta^*}(0,\gamma)}].$$

Since $\epsilon < \epsilon' \Rightarrow \mathcal{R}^{\alpha^* - \epsilon'}(\gamma + \delta_0) \subset \mathcal{R}^{\alpha^* - \epsilon}(\gamma + \delta_0)$,

$$\sup_{\epsilon' < \epsilon} \{\mathbb{1}_{0 \notin \mathcal{R}^{\alpha^* - \epsilon'}(\gamma + \delta_0)}\} = \mathbb{1}_{0 \notin \mathcal{R}^{\alpha^* - \epsilon}(\gamma + \delta_0)}$$

and (20) in Lemma 3 implies that, $\mu^{\alpha^*,\theta^*}$-a.e., $\lim_{\epsilon \to 0} \mathbb{1}_{0 \notin \mathcal{R}^{\alpha^* - \epsilon}(\gamma + \delta_0)} = 0$. Thus, according to Lebesgue's dominated convergence theorem and H4, $D_{21,n}$ vanishes asymptotically. For $D_{22,n}$, we apply H6$'$ and for $n$ sufficiently large,

$$D_{22,n} \leq \frac{\epsilon}{|\Lambda_n|} \sum_{x \in \gamma \cap \Lambda_n} \mathbb{1}_{\mathcal{R}^{\hat{\alpha}_n}(\gamma)}(x) |g((\gamma - \delta_x)_{-x})|.$$

According to the ergodic theorem, Proposition 2 and the stationarity of $\mu^{\alpha^*,\theta^*}$, we have, for $n$ sufficiently large,

$$D_{22,n} \leq 2\epsilon E_{\alpha^*,\theta^*}[|g(\gamma)| e^{-h^{\alpha^*,\theta^*}(0,\gamma)}],$$

which is less than $2\epsilon e^K E_{\alpha^*,\theta^*}[|g(\gamma)|]$. Therefore, $D_{22,n}$ also vanishes asymptotically and, as a consequence, $\lim_{n \to \infty} D_{2,n} = 0$.

Returning to (26), we have proved that both $D_{1,n}$ and $D_{2,n}$ vanish when $n \to \infty$. $\square$

## 5. Examples

In this section, we give three examples of Gibbs models which satisfy the assumptions of Section 4. The first is the classical model of hard spheres with pairwise step interaction and is hereditary. The second is the model of Delaunay tessellations with hardcore interaction on the size of the cells. It is a non-hereditary model and has been recently studied in [5]. The third is a model where each point interacts only with its $k$ nearest neighbors. A hardcore part forces these neighbors to be not too far away. These examples show that our assumptions are satisfied for a large class of models.

### 5.1. Hard sphere model with pairwise step interaction

Let $p$ be a positive integer and let $0 < r_1 < r_2 < \cdots < r_p$ be $p$ real numbers. For every $\alpha \in \mathcal{K}$ ($\mathcal{K}$ is a compact set in $\mathbb{R}^+ \setminus \{0\}$) and every $\theta = (\theta_1, \ldots, \theta_p)$ in $\Theta$ ($\Theta$ is a bounded open set in $\mathbb{R}^p$), the energy $H_\Lambda^{\alpha,\theta}(\gamma)$ is defined for every $\Lambda \in \mathcal{B}(\mathbb{R}^d)$ and $\gamma$ in $\mathcal{M}(\mathbb{R}^d)$ as

$$H_\Lambda^{\alpha,\theta}(\gamma) = \sum_{\substack{\{x,y\} \subset \gamma \\ \{x,y\} \cap \Lambda \neq \varnothing}} +\infty \mathbb{1}_{[0,1/\alpha]}(|x-y|) + \theta_1 \mathbb{1}_{]1/\alpha, 1/\alpha + r_1]}(|x-y|) + \cdots$$

$$+ \theta_p \mathbb{1}_{]1/\alpha + r_{p-1}, 1/\alpha + r_p]}(|x-y|). \tag{28}$$

*Remark 4.* We choose $\alpha$ in a compact set $\mathcal{K}$ to simplify some later calculations. In general, for the hard sphere model, the natural hardcore parameter is the radius of the hard sphere (see [6]). In our case, we take the inverse $\frac{1}{\alpha}$ so that S2 is fulfilled for $\alpha$.

Let us remark that the ranges of steps for the smooth part of the interaction depend on the parameter $\alpha$. This seems natural since the first step occurs after the hardcore radius $1/\alpha$ (otherwise, the model could be not identifiable) and, moreover, the lengths of the steps, $(r_1, \ldots, r_p)$, are fixed a priori.

This model is hereditary and the existence of Gibbs measures is well known (see [12], for example). Let us prove that it satisfies all the assumptions needed for our estimation procedure.

**Proposition 4.** *The family of energy functions $(H_\Lambda^{\alpha,\theta})$ for the hard sphere model satisfies assumptions* S1–S4*,* H2–H5 *and* H6′*.*

We first establish a lemma which is useful for checking S3.

**Lemma 5.** $\tilde{S}3$*, defined below, implies* S3*.*

$\tilde{S}3$. $\forall \alpha < \alpha^*$, $\exists R_\alpha > 0$, $\forall \gamma \in \mathcal{M}_\infty^{\alpha^*}(\mathbb{R}^d)$, $\exists (a_1, a_2, \ldots, a_m) \in B(0, R_\alpha)^m$ $(m \geq 1)$, $\exists \varepsilon > 0$, $\forall x_1 \in B(a_1, \varepsilon), \forall x_2 \in B(a_2, \varepsilon), \ldots, \forall x_m \in B(a_m, \varepsilon)$,

$$H_{B(0,R_\alpha)}^{\alpha^*,\theta^*}(\gamma_{B(0,R_\alpha)^c} + \delta_{x_1} + \cdots + \delta_{x_m}) < +\infty$$

*and*

$$H_{B(0,R_\alpha)}^{\alpha,\theta^*}(\gamma_{B(0,R_\alpha)^c} + \delta_{x_1} + \cdots + \delta_{x_m}) = +\infty.$$

**Proof.** Suppose that $(H_\Lambda^{\alpha,\theta^*})$ satisfies $\tilde{S}3$. Let $\alpha < \alpha^*$ and $\mu^{\alpha^*,\theta^*}$ be a Gibbs measure in $\mathcal{G}^{\alpha^*,\theta^*}$. We have

$$\mu^{\alpha^*,\theta^*}(H_{B(0,R_\alpha)}^{\alpha,\theta^*}(\gamma) = +\infty)$$

$$= \int \mathbb{1}_{\{H_{B(0,R_\alpha)}^{\alpha,\theta^*}(\gamma)=+\infty\}} \mu^{\alpha^*,\theta^*}(\mathrm{d}\gamma)$$

$$= \int \int \mathbb{1}_{\{H_{B(0,R_\alpha)}^{\alpha,\theta^*}(\gamma_{B(0,R_\alpha)^c}+\gamma')=+\infty\}} \frac{\mathrm{e}^{-H_{B(0,R_\alpha)}^{\alpha^*,\theta^*}(\gamma_{B(0,R_\alpha)^c}+\gamma')}}{Z_\Lambda^{\alpha^*,\theta^*}(\gamma_{B(0,R_\alpha)^c})} \pi_{B(0,R_\alpha)}(\mathrm{d}\gamma') \mu^{\alpha^*,\theta^*}(\mathrm{d}\gamma).$$

From $\tilde{S}3$, for $\mu^{\alpha^*,\theta^*}$ almost every $\gamma$ in $\mathcal{M}_\infty^{\alpha^*}(\mathbb{R}^d)$,

$$\int \mathbb{1}_{\{H_{B(0,R_\alpha)}^{\alpha,\theta^*}(\gamma_{B(0,R_\alpha)^c}+\gamma')=+\infty\}} \frac{\mathrm{e}^{-H_{B(0,R_\alpha)}^{\alpha^*,\theta^*}(\gamma_{B(0,R_\alpha)^c}+\gamma')}}{Z_\Lambda^{\alpha^*,\theta^*}(\gamma_{B(0,R_\alpha)^c})} \pi_{B(0,R_\alpha)}(\mathrm{d}\gamma') > 0.$$

Therefore, $\mu^{\alpha^*,\theta^*}(H_{B(0,R_\alpha)}^{\alpha,\theta^*}(\gamma) = +\infty) > 0$. The lemma is thus proved. $\square$

**Proof of Proposition 4.** S1 and S2 are obviously satisfied. The choice of $\frac{1}{\alpha}$ instead of $\alpha$ for the hardcore radius parameter is crucial here.

Choosing $R_\alpha = \frac{2}{\alpha}$, $r = \frac{1}{4}(\frac{1}{\alpha'} - \frac{1}{\alpha})$, $m = 2$, $a_1 = (0,0)$ and $a_2 = (\frac{1}{\alpha'} + 2r, 0)$, we check that $\tilde{S}3$ in Lemma 5 is satisfied, which implies S3.

S4 is satisfied because the number of terms in (28) is finite. So, given $\gamma$ such that $H_\Lambda^{\alpha',\theta}(\gamma) < +\infty$ for some $\alpha'$ in $\mathcal{K}$, it is easy to find $\alpha < \alpha'$ such that $H_\Lambda^{\alpha,\theta}(\gamma)$ remains finite.

Moreover, as the energy is clearly invariant by translation, H2 is obvious. Let us note that for every $\gamma \in \mathcal{M}_\infty^\alpha(\mathbb{R}^d)$,

$$h^{\alpha,\theta}(0,\gamma) = \sum_{x \in \gamma} +\infty \mathbb{1}_{[0,1/\alpha]}(|x|) + \theta_1 \mathbb{1}_{]1/\alpha,1/\alpha+r_1]}(|x|) + \cdots + \theta_p \mathbb{1}_{]1/\alpha+r_{p-1},1/\alpha+r_p]}(|x|). \qquad (29)$$

It is easy to see that $h^{\alpha,\theta}(0,\gamma)$ is either equal to infinity or uniformly bounded. Consequently, assumptions H3 and H4 are satisfied.

H5 is deduced from a classical identification property of the exponential family models (see, for example, assumption [ident] on page 244 of [2]).

It remains to prove H6'. For every $\gamma$ in $\mathcal{M}_\infty^{\alpha^*}(\mathbb{R}^d)$, we define

$$d(\gamma) = \min_{x \in \gamma}\left[\min\left(\left||x| - \left(\frac{1}{\alpha^*} + r_1\right)\right|, \left||x| - \left(\frac{1}{\alpha^*} + r_2\right)\right|, \ldots, \left||x| - \left(\frac{1}{\alpha^*} + r_p\right)\right|\right)\right].$$

$d(\gamma)$ is the minimal distance, over every $x$ in $\gamma$, between $|x|$ and a discontinuity point of the pairwise potential in (29). Since $\mu^{\alpha^*,\theta^*}$ is locally absolutely continuous with respect to the Poisson process, the random variable $d$ is $\mu^{\alpha^*,\theta^*}$-a.s. positive. Consequently, there exists a positive function $\varphi$ on $\mathbb{R}^+ \setminus \{0\}$ such that $\varphi$ is decreasing, $\lim_{x\to 0^+} \varphi(x) = +\infty$ and $\varphi(d(\cdot))$ is in $L^1(\mu^{\alpha^*,\theta^*})$. It suffices to choose, for instance, $\varphi = \frac{1}{\sqrt{F}}$, where $F$ is the repartition function of $d(\cdot)$ for the probability $\mu^{\alpha^*,\theta^*}$.

From the expression (29) of $h^{\alpha,\theta}(0,\gamma)$, there exists $K > 0$ such that, for every $\alpha \leq \alpha^*$ and $\gamma$ in $\mathcal{M}_\infty^{\alpha^*}(\mathbb{R}^d)$ which satisfy $0 \in \mathcal{R}^\alpha(\gamma + \delta_0)$ and $h^{\alpha,\theta}(0,\gamma) < +\infty$, we have, for every $(\theta,\theta')$ in $\Theta^2$,

$$|h^{\alpha,\theta}(0,\gamma) - h^{\alpha^*,\theta'}(0,\gamma)| \leq K \mathbb{1}_{[d(\gamma),+\infty[}\left(\frac{1}{\alpha} - \frac{1}{\alpha^*}\right) + K|\theta - \theta'|$$

$$\leq K \frac{\varphi(d(\gamma))}{\varphi(1/\alpha - 1/\alpha^*)} + K|\theta - \theta'|$$

$$\leq K \max(\varphi(d(\gamma)),1)\left[\frac{1}{\varphi(K'|\alpha - \alpha^*|)} + |\theta - \theta'|\right],$$

where $K'$ is a positive constant. This inequality implies H6$'$. Proposition 4 is thus proved. □

### 5.2. Hardcore Delaunay tessellations

In this section, we present the model of Delaunay tessellations studied in [5]. A hardcore interaction forces the cells of the tessellation to be not too small and not too large. Let us state some definitions concerning Delaunay tessellations. We assume that $d = 2$. Let $\gamma$ be in $\mathcal{M}(\mathbb{R}^d)$. A subset $T$ of three points in $\gamma$ is a *Delaunay triangle* if the open circumscribed ball $B(T)$ of $T$ does not contain any point of $\gamma$. $Del(\gamma)$ denotes the set of Delaunay triangles of $\gamma$.

Let $r > 0$ be a positive real number. For every $\alpha \in \mathcal{K}$ ($\mathcal{K}$ is a compact set in $\mathbb{R}^+ \setminus \{0\}$) and every $\theta \in \Theta$ ($\Theta$ is a bounded open set in $\mathbb{R}$), we define the energy $H_\Lambda^{\alpha,\theta}(\gamma)$, for every $\Lambda$ in $\mathcal{B}(\mathbb{R}^d)$ and $\gamma$ in $\mathcal{M}(\mathbb{R}^d)$, as

$$H_\Lambda^{\alpha,\theta}(\gamma) = \sum_{\substack{T \in Del(\gamma) \\ B(T) \cap \Lambda \neq \varnothing}} \varphi(T) \tag{30}$$

with

$$\varphi(T) = \begin{cases} +\infty, & \text{if } R(T) \geq \alpha, \\ +\infty, & \text{if } l(T) \leq r, \\ \theta Per(T), & \text{otherwise,} \end{cases} \tag{31}$$

where $R(T)$ is the radius of $B(T)$, $l(T)$ is the smallest length of the sides of $T$ and $Per(T)$ is the perimeter of $T$.

**Remark 5.** In this example, we choose to estimate the hardcore parameter $\alpha$ (as the maximum radius $R(T)$). It would also be possible to estimate the parameter $r$ (as the minimal length $l(t)$). Moreover, we believe that it is possible to estimate both of these parameters simultaneously, but this situation is not really dealt with in this paper and some further assumptions and proofs would be needed to deal with this case.

We have chosen the smooth part of the energy to be equal to $\theta Per(T)$ because it seems interesting to us. However, other smooth functions depending on the volume, the angles, etcetera are valuable and they could be managed in the same fashion.

**Proposition 5.** *The family of energy functions* $(H_\Lambda^{\alpha,\theta})$ *for the hardcore Delaunay tessellations model satisfies assumptions* S1–S4, H2–H5 *and* H6′ *if*

$$\forall \alpha \in \mathcal{K} \qquad \alpha > r. \tag{32}$$

**Proof.** S1, S2 and S4 are obviously satisfied, as in the example of Section 5.1. To prove S3, we use Lemma 5 and must prove that $\tilde{S}3$ is satisfied.

**Lemma 6.** $\tilde{S}3$ *is satisfied.*

**Proof.** Assumption (32) is crucial here. Let $\theta$ be in $\Theta$ and $\alpha < \alpha^*$. Let $\tilde{T} = \{a_1, a_2, a_3\}$ be an equilateral triangle in $\mathbb{R}^d$, where the radius $R(\tilde{T})$ of the circumscribed ball $B(\tilde{T})$ satisfies $\alpha < R(\tilde{T}) < \alpha^*$. Moreover, let $R_\alpha > 0$ be such that $B(\tilde{T}) \subset B(0, R_\alpha - r - \alpha^*)$. Now, for every $\gamma$ in $\mathcal{M}_\infty^{\alpha^*}(\mathbb{R}^d)$, we must add some points $a_4, a_5, \ldots, a_m$ to the points $a_1$, $a_2$, $a_3$ to obtain the expected configuration needed in assumption $\tilde{S}3$.

We construct these points $a_4, a_5, \ldots, a_m$ recursively. This construction is inspired by the one given in [5], page 148. Let us start with $\{a_4, a_5, \ldots, a_{m'}\} = \gamma_{B(0,R_\alpha) \setminus B(0,R_\alpha - \alpha^*)}$ and put $\hat{\gamma} = \gamma_{B(0,R_\alpha)^c} + \delta_{a_1} + \delta_{a_2} + \delta_{a_3} + \delta_{a_4} + \cdots + \delta_{a_{m'}}$. We now test whether $\hat{\gamma}$ satisfies the assumption

$$\forall x \in B(0, R_\alpha) \qquad \hat{\gamma}(\bar{B}(x, \alpha^*)) > 0, \tag{A}$$

where $\bar{B}(x, \alpha^*)$ denotes the closed ball centered at $x$ with radius $\alpha^*$. If this is the case, we stop the recursive process, put $m = m'$ and $a_1, a_2, \ldots, a_m$ are obtained. If (A) is false, we choose an arbitrary point $a_{m'+1}$ in $B(0, R_\alpha)$ such that $\hat{\gamma}(\bar{B}(a_{m'+1}, \alpha^*)) = 0$ and we put $\hat{\gamma} := \hat{\gamma} + \delta_{a_{m'+1}}$. We then test (A) for this new $\hat{\gamma}$. If (A) is true, we stop here and $m = m' + 1$; otherwise, we choose another point $\hat{a}_{m'+2}$ as above and put $\hat{\gamma} := \hat{\gamma} + \delta_{a_{m'+2}}$. We continue like this until the process stops, which is always the case since $B(0, R_\alpha)$ is a bounded set and, from (32), the minimal distance between two points in $\hat{\gamma}$ is larger than $r$. Finally, it is relatively easy to see that

$$H_{B(0,R_\alpha)}^{\alpha^*,\theta^*}(\gamma_{B(0,R_\alpha)^c} + \delta_{a_1} + \cdots + \delta_{a_m}) < +\infty$$

and

$$H_{B(0,R_\alpha)}^{\alpha,\theta^*}(\gamma_{B(0,R_\alpha)^c} + \delta_{a_1} + \cdots + \delta_{a_m}) = +\infty.$$

Both of these properties remain true for small local perturbations of $a_1, a_2, \ldots, a_m$. So, for $\varepsilon$ sufficiently small, $\tilde{S}3$ is proved. □

H2–H5 and H6′ still have to be proven. The energy is clearly invariant under translation, so H2 is obviously satisfied. Moreover, for every $\gamma \in \mathcal{M}_\infty^\alpha(\mathbb{R}^d)$,

$$h^{\alpha,\theta}(0,\gamma) = \sum_{\substack{T \in Del(\gamma+\delta_0) \\ 0 \in T}} \varphi(T) - \sum_{\substack{T \in Del(\gamma) \\ 0 \in B(T)}} \varphi(T), \tag{33}$$

so $h^{\alpha,\theta}(0,\gamma)$ is either equal to infinity or uniformly bounded (see [5], page 137). Consequently, assumptions H3 and H4 are satisfied. Contrary to the example of Section 5.1, assumption H6′ is obviously satisfied here. H5 still has to be proven. With a proof similar to the one given in Lemma 6 above, we show that $\forall \theta \in \Theta \setminus \{\theta^*\}, \exists R_{\alpha^*} > 0, \forall \gamma \in \mathcal{M}_\infty^{\alpha^*}(\mathbb{R}^d)$, $\exists (a_1, a_2, \ldots, a_m) \in B(0, R_{\alpha^*})^m, \exists \varepsilon > 0, \forall x_1 \in B(a_1,\varepsilon), \forall x_2 \in B(a_2,\varepsilon), \ldots, \forall x_m \in B(a_m,\varepsilon)$

$$h^{\alpha^*,\theta^*}(0, \gamma_{B(0,R_{\alpha^*})^c} + \delta_{x_1} + \cdots + \delta_{x_m}) \neq h^{\alpha^*,\theta}(0, \gamma_{B(0,R_{\alpha^*})^c} + \delta_{x_1} + \cdots + \delta_{x_m}). \tag{34}$$

With an argument similar to that used for Lemma 5, it can be proven that (34) implies H5. Proposition 5 is thus proved. □

### 5.3. A forced clustering $k$-nearest-neighbors model

In this section, we introduce a model where any point $x$ interacts with only its $k$ nearest neighbors (see [1], for example). We propose including a hardcore interaction which forces these $k$ nearest neighbors to be at a distance less than $\alpha$ from $x$. Let us describe the interaction precisely.

Let $k \geq 1$ be a fixed integer and let $\varphi$ be a bounded function from $\mathbb{R}^+$ to $\mathbb{R}$ which is non-null in a neighborhood of 0.

For every $\alpha \in \mathcal{K}$ ($\mathcal{K}$ is a compact set in $\mathbb{R}^+ \setminus \{0\}$) and every $\theta \in \Theta$ ($\Theta$ is a bounded open set in $\mathbb{R}$), we define the energy $H_\Lambda^{\alpha,\theta}(\gamma)$, for every $\Lambda \in \mathcal{B}(\mathbb{R}^d)$ and every $\gamma$ in $\mathcal{M}(\mathbb{R}^d)$, as

$$H_\Lambda^{\alpha,\theta}(\gamma) = \sum_{x \in \gamma_{\Lambda^\alpha}} \left[ +\infty \mathbb{1}_{[0,k]}(\gamma(\bar{B}(x,\alpha))) + \sum_{y \in \mathcal{N}^k(x,\gamma)} \theta\varphi(|x-y|) \right], \tag{35}$$

where $\Lambda^\alpha$ denotes the set $\bigcup_{x \in \Lambda} B(x,\alpha)$, $\mathcal{N}^k(x,\gamma)$ is the set containing the $k$ nearest neighbors of $x$ in $\gamma$ and $\bar{B}(x,\alpha)$ is the closed ball centered at $x$ with radius $\alpha$.

The hardcore interaction implies that the locally finite energy configurations do not contain clusters with less than $k+1$ points. This model therefore naturally forces clusters and is a non-hereditary model.

In the literature, no proof for the existence of Gibbs measures for these kinds of interactions is available, but we claim that the proof is essentially the same as that given for Example 2 in [5].

**Proposition 6.** *The family of energy functions $(H_\Lambda^{\alpha,\theta})$ for the forced clustering k-nearest-neighbors model satisfies assumptions* S1–S4*,* H2–H5 *and* H6′*.*

**Proof.** S1, S2 and S4 are satisfied as in both the examples of Sections 5.1 and 5.2. S3 is satisfied by checking S̃3 in Lemma 5, which can be easily done. For the other assumptions, note that for every $\gamma \in \mathcal{M}_\infty^\alpha(\mathbb{R}^d)$,

$$h^{\alpha,\theta}(0,\gamma) = +\infty \mathbb{1}_{[0,k]}(\gamma(\bar{B}(0,\alpha))) + \sum_{x \in \mathcal{N}^k(0,\gamma+\delta_0)} \theta\varphi(|x|) \\ + \sum_{\substack{x \in \gamma \\ 0 \in \mathcal{N}^k(x,\gamma+\delta_0)}} \theta\varphi(|x|) - \theta\varphi(|\chi^k(x,\gamma)|), \qquad (36)$$

where $\chi^k(x,\gamma)$ is the farthest neighbor of $x$ in $\mathcal{N}^k(x,\gamma)$.

In (36), the number of terms in the first sum is obviously equal to $k$. In the second sum, the number of terms is uniformly bounded (for instance, by $5k$ when $d=2$). Consequently, as in both the examples of Sections 5.1 and 5.2, $h^{\alpha,\theta}(0,\gamma)$ is either equal to infinity or uniformly bounded. Therefore, assumptions H3, H4 and H6′ are satisfied and H5 is proved by an argument similar to that used in the example of Section 5.2. □